\documentclass[a4paper,12pt]{article}
\usepackage{amsmath}
\usepackage{amssymb}
\usepackage{tabularx}
\usepackage{enumerate}
\usepackage{graphicx}
\usepackage{subfigure}
\usepackage{epsfig}
\usepackage{graphics}
\usepackage{cases}
\usepackage[compress]{cite}
\usepackage{txfonts}
\usepackage{geometry}
 \usepackage{epstopdf}

\newtheorem{thm}{Theorem}[section]
\newtheorem{lem}[thm]{Lemma}

\newtheorem{exm}{Example}[section]

\newtheorem{rmk}{Remark}[section]

\newtheorem{pppp}{Proof}

\newcommand{\qed}{\hspace{1em}\mbox{\raisebox{0.65ex}{\fbox{}}}}

\numberwithin{equation}{section}

\newcommand{\be}{\begin{equation}}
\newcommand{\ee}{\end{equation}}
\newcommand\bes{\begin{eqnarray}} \newcommand\ees{\end{eqnarray}}
\newcommand{\bess}{\begin{eqnarray*}}
\newcommand{\eess}{\end{eqnarray*}}

\newcommand{\bpf}{{\bf Proof:\ \ }}
\newcommand{\epf}{\mbox{}\hfill $\Box$}

\begin{document}

\thispagestyle{empty}

\title{Effects of depth and evolving rate on phytoplankton growth in a periodically evolving environment\thanks{The work is partially supported by the NNSF of China (Grant No. 11771381) and Research Foundation of Young Teachers of Hexi University (QN2018013).}}

\date{\empty}

\author{Liqiong Pu$^{a, b}$ and Zhigui Lin$^a\thanks{Corresponding author. Email: zglin@yzu.edu.cn (Z. Lin).}$\\
{\small $^a$ School of Mathematical Science, Yangzhou University, Yangzhou 225002, China}\\
{\small $^b$ School of Mathematics and Statistics, Hexi University, Zhangye 734000, China.}
}

 \maketitle
\begin{quote}
\noindent
{\bf Abstract.} { 
 In this paper, we consider a single phytoplankton species which relies on the light for maintaining the metabolism of life in a periodically evolving environment, where the light intensity and the death rate depend on the water column depth triggered by seasonal variation. Based on the basic reproduction number $\mathcal{R}_0$, a threshold type result on  the dynamics of the model is established. Especially, various features of $\mathcal{R}_0$ with respect to the vertical turbulent diffusion rate, the buoyant or sinking rate, and the evolving rate of water column depth are derived. Our theoretical results and numerical simulations show that big evolving rate, vertical diffusion rate and water column depth all have an adverse effect on survival of phytoplankton.
}

\noindent {\it MSC:} 35K57, 37B55, 92D30

\medskip
\noindent {\it Keywords:}  Phytoplankton growth;  Evolving environment; Basic reproduction number; Threshold dynamics.
\end{quote}

\section{Introduction}

Phytoplankton populations include a diverse representation of both microscopic eukaryote algae and prokaryote cyanobacteria. They are major oxygen producers and also represent the initial food source for numerous food webs in freshwater, estuarine, and marine waters \cite{MHG}. Since they transport significant amounts of atmospheric carbon dioxide into the deep oceans, they play a crucial role in the climate dynamics. Its growth and reproduction are not only influenced by its own biological characteristics, but also by substantial environmental factors such that light, water flow, temperature, salinity, nutrient, stress and disturbance, et al.

Influence mechanism of the
above factors on phytoplankton growth is an essential and important subject in mathematical ecology.
Edwards et al. \cite{ET} proposed that light limitation diminishes the temperature sensitivity of bulk phytoplankton growth, even though community structure will be temperature-sensitive. To describe the interactions between photosynthesis and light for phytoplankton,
Jassby and Platt \cite{JP} gave a mathematical formulation, which is most consistent with the data. A mathematical model
was proposed by Rhee and Gotham \cite{wy} to characterize the growth of two types of phytoplankton, non-nitrogen-fixing and nitrogen-fixing, which both require light in order to grow. Considering the effect of flow velocity,
Li et al. \cite{lzz} showed that phytoplankton biomass and spatial distribution largely depend on the flow condition by the field observation, and turbulent flow has the inhibition effect on phytoplankton biomass, but less impact on composition by enclosure experiments.
Chellappa et al. \cite{CN} studied a phytoplankton growth with stress and disturbance, which is stimulated by variable hydrology, and indicated that stress and disturbance factors determine phytoplankton diversity, composition and chlorophyll levels.

Temperature, as a physical factor, is of supreme significance to affect the growth of phytoplankton. The earlier work was discussed by Canale and Vogel \cite{CV}.
Later, Rhee and Gotham \cite{RG} investigated the simultaneous effects of nutrient limitation and suboptimal temperature on growth and nutrient uptake kinetics by continuous-culture techniques. Based on patterns observed in nature, Agawin et al. \cite{ADA} studied the combined effects of both temperature and nutrient concentrations, and tested the relative importance of phytoplankton in warm and nutrient-poor waters. Considering phytoplankton living in freshwater and brackish water,
Lionard et al. \cite{LMGV} concluded that increasing or decreasing of salinity brings negative impart on the growth of phytoplankton.

In this paper, we consider the self-shading model with unlimited nutritional environment proposed by Shigesada and Okubo \cite{SO}, the governing equation is
$$u_{t}=Du_{xx}-\alpha u_x+F(u)-\lambda u,$$
where $F(u)=\int_0^{u}f(\frac{I_0}{I_m}e^{-x})dx$, the function $f$ describes the dependence of production rate, and depends on the incident light intensity $I_0$ and a characteristic light intensity $I_m$. $\lambda$ is the ratio of mortality to maximum growth rate.
The model has been receiving extensive investigations, see for example, Ebert et al. \cite{EAT}, Kolokolnikov, Ou and Yuan \cite{ky}, Du and Hsu \cite{DH}, and Hsu and Lou \cite{HL} and the references therein.

In order to understand the evolution of a single phytoplankton species which depends merely on light for its metabolism in a eutrophic and vertical water bodies,
recently, Peng and Zhao \cite{pzh} developed the model and studied a nonlocal reaction-diffusion-advection model of a single phytoplankton species
\begin{eqnarray}
\left\{
\begin{array}{lll}
u_{t}=Du_{xx}-\alpha u_x+u[g(I(x,t))-d(x,t)],\; &\,0<x<L_0,\ t>0, \\[2mm]
Du_x(x,t)-\alpha u(x,t)=0,\; &\, x=0,L_0,\ t>0,\\[2mm]
u(x,0)=u_0(x)\geq,\not\equiv0,\; &\, 0<x<L_0,\ t>0,
\end{array} \right.
\label{a01}
\end{eqnarray}
where light intensity $I(x,t)=I_0e^{-k_0x-k_1\int_{0}^{x}u(s,t)ds}$, $I_0=I_0(t)\geq 0, \not\equiv 0$, $I_0$ is $T-$perodic in $t$ for some $T>0$ and denotes the incident light intensity, $k_0$ and $k_1$ are positive constants which express the background turbidity and the absorption coefficient of the phytoplankton species respectively. Specific growth rate $g(I)$ is a function of the light intensity $I(x,t)$ and $g$ is strictly increasing with $g(0)=0$. The death rate $d(x,t)$ is a $T-$periodic in $t$ and satisfies that either $d_x(x,t)\geq0,\not\equiv0$ on $[0,L_0]\times[0,T]$ and $k\geq0$, or $d(x,t)\equiv d(t)$ and $k_0>0$. $x$ stands for the depth from $0$ (the top) to $L_0$ (the bottom)
in a column of water. They investigated persistence and extinction of phytoplankton species in accordance with the basic reproduction number $\mathcal{R}_0$, which depends
on the water column depth $L$, the vertical turbulent diffusion rate $D$ and the settling (or buoyant) rate $\alpha>0( $or $<0 )$.

\begin{figure}[ht]
\centering
\subfigure[]{ {
\includegraphics[width=0.30\textwidth]{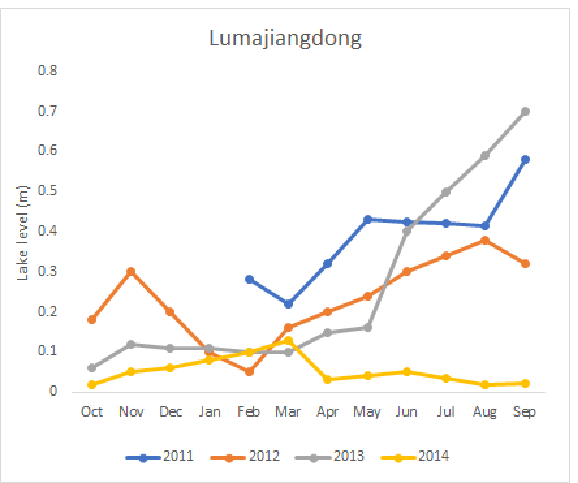}
} }
\subfigure[]{ {
\includegraphics[width=0.30\textwidth]{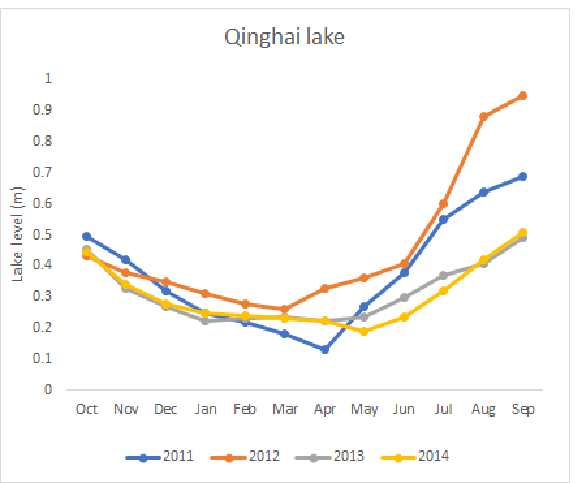}
} }
\caption{Figures (a) and (b) describe level changes of Lumajiangdong and Qinghai lake in a hydrological year (October to September) from 2011 to 2014, respectively. These data are from \cite{LY} (minor modification).
}
\label{tu01}
\end{figure}

As is well-known, lake depth changes periodically, and the periodic annual changes mainly depend on the replenishment of lake water. Lakes supply comes from precipitation, such as Poyang lake, which is China's largest fresh water lake, has the highest water level in rainy season and the lowest in dry season \cite{YLL}. Plateau lakes are mainly replenished by melting water from snow and ice, with the highest water level in summer and the lowest in winter. Take lumajiangdong and qinghai lake for examples. Trends in lake level of lakes across  the Tibetan Plateau in a hydrological year are given in \cite{LY}, it can be easily seen from Fig. \ref{tu01} that the lakes have striking characteristics, their lake level in summer is significantly higher than that in winter.

Moreover, it is widely understood that evolving of environment plays a critical role in species' migration \cite{C1, CGM2, CM}.
The rivers, lakes and seas, where phytoplankton species habitats, are evolving, therefore studying the dispersal of phytoplankton in an evolving environment is essential for understanding phytoplankton growth in biological systems.

As in \cite{C1, tql1}, let $(0,L(t))$ be a bounded evolving interval at time $t>0$ with its shifting depth $L(t)$. For any point $x(t) \in (0,L(t))$, we assume that $u(x(t),t)$ is the population density of the phytoplankton species at position $x(t)$ and time $t> 0$.

Let $x(t)=a(t)$ and $x(t)=b(t)( a(t)\leq b(t) )$ denote arbitrary two depth within the water column where $x(t)$ varies from $0$ ( the top ) to $L(t)$ ( the bottom ), obviously, $(a(t),b(t))\subset(0,L(t))$ be a bounded evolving interval at time $t>0$ with its moving top $x(t)=a(t)$ and bottom $x(t)=b(t)$.
According to the law of mass conservation, one sees immediately that
\bes
\begin{array}{lll}
\frac{d}{dt}\int_{a(t)}^{b(t)}u(x(t),t)dx
&=&Du_x(b(t),t)-Du_x(a(t),t)\\[2mm]
& &-\alpha\int_{a(t)}^{b(t)}u_x(x(t),t)dx+\int_{a(t)}^{b(t)}f(u(x(t),t))dx\\[2mm]
&=&\int_{a(t)}^{b(t)}(Du_{xx}-\alpha u_x+f(u))dx,
\end{array}
\label{a02}
\ees
where $\frac{d}{dt}\int_{a(t)}^{b(t)}u(x(t),t)dx$ represents the change rate of phytoplankton population density in the interval $(a(t),b(t))$, the term $Du_x(b(t),t)-Du_x(a(t),t)$ denotes the net flow through the boundary (flows in from boundary $x=a(t)$ and out of boundary $x=b(t)$), $-\alpha\int_{a(t)}^{b(t)}u_x(x(t),t)dx$ is the net increase of the advection term, and $\int_{a(t)}^{b(t)}f(u(x(t),t))dx$ represents the net increase of the reaction term in the interval $(a(t),b(t))$.

Equation \eqref{a02} can be simplified if we suppose that the region evolution is isotropic and uniform.
That is, the evolution of the interval takes place at the same proportion as time elapses.
Mathematically, $x(t)$ can be described as follows
\bes
x(t)=\rho(t)y, &y\geq 0,
\label{a03}
\ees
where the positive continuous function $\rho(t)$ is called evolving rate subject to $\rho(0)=1$. Let $L_0=L(0)$, we have evolving interval  $(0,L(t))=(0,\rho(t)L_0)$. Set $a(t)=\rho(t)y_1$ and $b(t)=\rho(t)y_2$, therefore $y_1,y_2\in (0,L_0)$. Furthermore, if $\rho(t)=\rho(t+T)$ for some $T>0$, the habitat is periodically evolving, while if $\dot \rho(t) \ge0$, we say it is growing and if $\dot \rho(t) \le0$, it is called shrinking habitat.

Noticing that
$u(x(t),t)=u(\rho(t)y,t)$, we then define $u(x(t),t)=v(y,t)$.
From the left side of \eqref{a02}, we get
\bes
\begin{array}{llll}
\frac{d}{dt}\int_{a(t)}^{b(t)}u(x(t),t)dx&=&\frac{d}{dt}\int_{\rho(t)y_1}^{\rho(t)y_2}u(\rho(t)y,t)d(\rho(t)y)\\[2mm]
&=&\frac{d}{dt}\int_{y_1}^{y_2}v(y,t)\rho(t)dy\\[2mm]
&=&\int_{y_1}^{y_2}(v_t(y,t)\rho(t)+v(y,t)\dot\rho(t))dy.
\end{array}
\label{a04}
\ees
Using \eqref{a03} yields
\bess
&&u_{xx}=\frac{1}{\rho^2(t)}v_{yy},\quad
u_x=\frac{1}{\rho(t)}v_y,
\eess
then the right side of \eqref{a02} becomes
\bes
\begin{array}{llll}
\int_{a(t)}^{b(t)}(Du_{xx}-\alpha u_x+f(u) )dx=\int_{y_1}^{y_2}(\frac{D}{\rho^2(t)}v_{yy}-\frac{\alpha}{\rho(t)}v_y+f(v(y,t)))\rho(t)dy.
\end{array}
\label{a05}
\ees
We conclude from \eqref{a02}, \eqref{a04} and \eqref{a05} that
\bes
\int_{y_1}^{y_2}(v_t+\frac{\dot{\rho}(t)}{\rho(t)}v)dy=\int_{y_1}^{y_2}(\frac{D}{\rho^2(t)}v_{yy}
-\frac{\alpha}{\rho(t)}v_y + f(v,y,t))dy,  & \, t>0.
\label{a06}
\ees
According to arbitrariness of $a(t)$ and $b(t)$, \eqref{a06} holds for any $y_1, y_2\in (0,L_0)$. On account of
its boundary condition, we derive the following problem in a fixed interval
\bes\left\{
\begin{array}{llll}
v_{t}=\frac{D}{\rho^2(t)} v_{yy}-\frac{\alpha}{\rho(t)}v_y+[g(I(\rho(t)y,t))-d(\rho(t)y,t)-\frac{\dot{\rho}(t)}{\rho(t)}]v, \;& 0<y<L_0, \; t>0,\\
Dv_y-\alpha\rho(t)v(y,t)=0,\;& y=0,L_0, \;  t>0
\end{array}\right.
\label{a07}
\ees
with the initial condition
\begin{equation}
v(y,0)=u_0(y)\geq 0,\not\equiv 0, \quad 0<y<L_0.
\label{a08}
\end{equation}

For later application, we also consider problem \eqref{a07} with the periodic condition
\begin{equation}
v(y,0)=v(y,T), \quad 0<y<L_0.
\label{a09}
\end{equation}
In order to study the impact of $D$ and $\alpha$, we would better make a transformation. Let $z(y,t)=e^{-\frac{\alpha}{D}\rho(t)y}v(y,t)$, \eqref{a07} with \eqref{a08} turns to
{\small
\bes\left\{
\begin{array}{llllll}
z_{t}=\frac{D}{\rho^2(t)} z_{yy}+\frac{\alpha}{\rho(t)}z_y+[g(I(\rho(t)y,t))
-d(\rho(t)y,t)-\frac{\dot{\rho}(t)}{\rho(t)}-\frac{\alpha y\dot{\rho}(t)}{D}]z, \;& 0<y<L_0, \; t>0,\\
z_y(y,t)=0,\;& y=0,L_0, \;  t>0,\\
z(y,0)=u_0(y)\geq 0,\not\equiv 0, \; & 0<y<L_0.
\end{array}\right.
\label{a10}
\ees
}

The present paper is build up as follows. In Section 2, we define the basic reproduction number based on the principal eigenvalue, a comparison principle for the basic reproduction number is established. Moreover, various features of $\mathcal{R}_0$ with respect to the vertical turbulent diffusion rate and the water column depth are derived. Section 3 is devoted to the asymptotic behavior of solution. The paper ends with some simulations and ecological significance for our analytical findings.

\section{The basic reproduction number and comparison principle}

In this section, we define the basic reproduction number ${\mathcal{R}}_0$, and then investigate its threshold-type dynamics for the system \eqref{a07}, \eqref{a08} in terms of ${\mathcal{R}}_0$. To do so, we first consider the corresponding periodic problem.

Obviously, $v^*\equiv0$ is an equilibrium solution of the periodic problem \eqref{a07}, \eqref{a09}, a routine computation gives rise to the corresponding linearized system of problem \eqref{a07} at $v^*$,
\begin{eqnarray}
\left\{
\begin{array}{lll}
w_{t}=\frac{D}{\rho^2(t)} w_{yy}-\frac{\alpha}{\rho(t)}w_y+g(I_0e^{-k_0\rho(t)y})w
-(d(\rho(t)y,t)+\frac{\dot{\rho}(t)}{\rho(t)})w, \;& 0<y<L_0, \; t>0,\\
Dw_y-\alpha\rho(t)w(y,t)=0,\;& y=0,L_0, \;  t>0.
\end{array} \right.
\label{b01}
\end{eqnarray}
Let $\eta(t,s)$ be the evolution operator of the problem
\begin{eqnarray}
\left\{
\begin{array}{lll}
w_{t}=\frac{D}{\rho^2(t)} w_{yy}-\frac{\alpha}{\rho(t)}w_y-(d(\rho(t)y,t)+\frac{\dot{\rho}(t)}{\rho(t)})w, \;& 0<y<L_0, \; t>0,\\
Dw_y-\alpha\rho(t)w(y,t)=0,\;& y=0,L_0, \;  t>0,
\end{array} \right.
\label{b02}
\end{eqnarray}
because of $d(\rho(t)y,t)>0$ for $(y,t)\in[0,L_0]\times[0,T]$, $\int_{s}^{t}\frac{\dot{\rho}(z)}{\rho(z)}dz$ is bounded for $\forall \,t\geq s$. On account of the standard
semigroup theory, it has become obviously that there exist positive constants $K$ and $c_0$ such that
$$\|\eta(t,s)\|\leq Ke^{-c_0(t-s)},\ \ \  \forall \, t\geq s,\,\, t, s\in \mathbb{R}.$$
Let $C_T$ be the ordered Banach space consisting of all  $T-$ periodic
and continuous function from $\mathbb{R}$ to $C([0,L_0],\mathbb{R})$ with the maximum norm $\|\cdot\|$ and the positive cone $C_T^+:=\{\mu\in C_T:\sigma(y)\geq0, \forall \ t\in \mathbb{R}, y\in[0,L_0]\}$. The notation $\sigma(y,t):=\sigma(t)(y)$ will be adopted for any given $\sigma\in C_T.$ According to \cite{ZXQ}, a linear operator $A$ on $C_T$ is introduced
$$ A(\sigma)(t):=\int_{0}^{\infty}\eta(t,t-s)g(I_0(t-a)e^{-k_0\cdot})\sigma(\cdot,t-s)ds. $$
From the hypothesis of $g$ and $I_0$, it is quite clear that $A$ is positive, continuous and compact on $C_T$. We define the spectral radius of $A$
$$ {\mathcal{R}}_0=\rho(A) $$
as the basic reproduction number for periodic system \eqref{a07},\eqref{a09}.
In order to describe the basic reproduction number more clearly, the following two eigenvalue problems are given
\begin{eqnarray}
\left\{
\begin{array}{lll}
\phi_{t}=\frac{D}{\rho^2(t)} \phi_{yy}-\frac{\alpha}{\rho(t)}\phi_y+[\frac{g(I_0(t)e^{-k_0\rho(t)y})}{\mu}
-d(\rho(t)y,t)-\frac{\dot{\rho}(t)}{\rho(t)}]\phi,\;& 0<y<L_0, \; 0<t\leq T,\\
D\phi_y-\alpha\rho(t)\phi(y,t)=0,\;& y=0,L_0, \;  0<t\leq T,\\
\phi(y,0)=\phi(y,T),\;& 0<y<L_0,
\end{array} \right.
\label{b03}
\end{eqnarray}
and
\begin{eqnarray}
\left\{
\begin{array}{lll}
\psi_{t}=\frac{D}{\rho^2(t)} \psi_{yy}-\frac{\alpha}{\rho(t)}\psi_y+[g(I_0(t)e^{-k_0\rho(t)y}&\\
\quad \quad -d(\rho(t)y,t))-\frac{\dot{\rho}(t)}{\rho(t)}]\psi+\lambda\psi,\;& 0<y<L_0, \; 0<t\leq T,\\
D\psi_y-\alpha\rho(t)\psi(y,t)=0,\;& y=0,L_0, \;  0<t\leq T,\\
\psi(y,0)=\psi(y,T),\;& 0<y<L_0,
\end{array} \right.
\label{b04}
\end{eqnarray}
where $(\mu_0,\phi)$ and $(\lambda_0,\psi)$ are the principal eigen-pairs of eigenvalue problems \eqref{b03} and \eqref{b04}, respectively. With above definitions, we have the following  result, see similar results in \cite{pzh}.
\begin{lem}
$(i)$  ${\mathcal{R}}_0=\mu_0$.
$(ii)\  \rm{sign} (1-\mathcal{R}_0)=\rm{sign}\lambda_0$.
\end{lem}

Further, based on \eqref{a10}, $(\mathcal{R}_0,\varphi(y,t))$ with $\varphi(y,t)=e^{-\frac{\alpha}{D}\rho(t)y}\phi(y,t)$ is the principal eigen-pair of the following
periodic-parabolic eigenvalue problem
\bes\left\{
\begin{array}{llllll}
\varphi_{t}=\frac{D}{\rho^2(t)} \varphi_{yy}+\frac{\alpha}{\rho(t)}\varphi_y+[\frac{g(I_0(t)e^{-k_0\rho(t)y})}{\mathcal{R}_0}\\[2mm]
\qquad -d(\rho(t)y,t)-\frac{\dot{\rho}(t)}{\rho(t)}-\frac{\alpha y\dot{\rho}(t)}{D}]\varphi, \;& 0<y<L_0, \; 0<t\leq T,\\
\varphi_y(y,t)=0,\;& y=0,L_0, \;  0<t\leq T,\\
\varphi(y,0)=\varphi(y,T),\;& 0<y<L_0,
\end{array}\right.
\label{b05}
\ees
which can be written as the following
 \bes\left\{
\begin{array}{llllll}
(e^{\frac{\alpha}{D}\rho(t)y}\varphi)_{t}=\frac{D}{\rho^2(t)}(e^{\frac{\alpha}{D}\rho(t)y} \varphi_{y})_y+
[\frac{g(I_0(t)e^{-k_0\rho(t)y})}{\mathcal{R}_0}\\[2mm]
\qquad \qquad  \ \  \ -d(\rho(t)y,t)-\frac{\dot{\rho}(t)}{\rho(t)}]e^{\frac{\alpha}{D}\rho(t)y}\varphi, \;& 0<y<L_0, \; 0<t\leq T,\\
\varphi_y(y,t)=0,\;& y=0,L_0, \;  0<t\leq T,\\
\varphi(y,0)=\varphi(y,T),\;&  0<y<L_0.
\end{array}\right.
\label{b06}
\ees
\begin{lem}(Comparison principle)

$(i)$ If there is a nonnegative nontrivial function $\tilde{z}\in C^{2,1}([0,L_0]\times[0,T])$, such that
\bes\left\{
\begin{array}{llllll}
(e^{\frac{\alpha}{D}\rho(t)y}\tilde{z})_{t}\geq\frac{D}{\rho^2(t)}(e^{\frac{\alpha}{D}\rho(t)y} \tilde{z}_{y})_y+
[\frac{g(I_0(t)e^{-k_0\rho(t)y})}{\mathcal{R}_1}\\[2mm]
\qquad \qquad  \ \  \ \ -d(\rho(t)y,t)-\frac{\dot{\rho}(t)}{\rho(t)}]e^{\frac{\alpha}{D}\rho(t)y}\tilde{z}, \;& 0<y<L_0, \;0<t\leq T,\\
\tilde{z}_y(0,t)\leq0,\;\tilde{z}_y(L_0,t)\geq0\; &  0<t\leq T, \\
\tilde{z}(y,0)\geq\tilde{z}(y,T),\;&  0<y<L_0,
\end{array}\right.
\label{b07}
\ees
for some real number $\mathcal{R}_1$, then $\mathcal{R}_0\leq\mathcal{R}_1$ and the equality holds only when $(\mathcal{R}_0,\tilde{z})$ is the principal eigen-pair of eigenvalue problem \eqref{b05}.

$(ii)$ If there is a nonnegative nontrivial function $\hat{z}\in C^{2,1}([0,L_0]\times[0,T])$, such that
\bes\left\{
\begin{array}{llllll}
(e^{\frac{\alpha}{D}\rho(t)y}\hat{z})_{t}\leq\frac{D}{\rho^2(t)}(e^{\frac{\alpha}{D}\rho(t)y} \hat{z}_{y})_y+
[\frac{g(I_0(t)e^{-k_0\rho(t)y})}{\mathcal{R}_2}\\[2mm]
\qquad \qquad  \ \  \ -d(\rho(t)y,t)-\frac{\dot{\rho}(t)}{\rho(t)}]e^{\frac{\alpha}{D}\rho(t)y}\hat{z}, \;& 0<y<L_0, \; 0<t\leq T,\\
\hat{z}_y(0,t)\geq0,\;\hat{z}_y(L_0,t)\leq0\; &  0<t\leq T, \\
\hat{z}(y,0)\leq\hat{z}(y,T),\;&  0<y<L_0,
\end{array}\right.
\label{b08}
\ees
for some real number $\mathcal{R}_2$, then $\mathcal{R}_0\geq\mathcal{R}_2$ and the equality holds only when $(\mathcal{R}_0,\hat{z})$ is the principal eigen-pair of eigenvalue problem \eqref{b05}.
\end{lem}
 \bpf
We only verify the assertion $(i)$, because $(ii)$ can be proved in the same manner.
From \cite{hp}, one knows that for the adjoint problem of \eqref{b06}
\bes\left\{
\begin{array}{llllll}
-e^{\frac{\alpha}{D}\rho(t)y}w_{t}=\frac{D}{\rho^2(t)}(e^{\frac{\alpha}{D}\rho(t)y} w_{y})_y+
[\frac{g(I_0(t)e^{-k_0\rho(t)y})}{\mathcal{R}_0}\\[2mm]
\qquad \qquad \; \; \ \  \ -d(\rho(t)y,t)-\frac{\dot{\rho}(t)}{\rho(t)}]e^{\frac{\alpha}{D}\rho(t)y}w, \;& 0<y<L_0, \; 0<t\leq T,\\
w_y(y,t)=0,\;& y=0,L_0, \;  0<t\leq T,\\
w(y,0)=w(y,T),\;&  0<y<L_0,
\end{array}\right.
\label{b09}
\ees
which is equivalent to
\bes\left\{
\begin{array}{llllll}
-w_{t}=\frac{D}{\rho^2(t)}w_{yy}+\frac{\alpha}{\rho(t)}w_y+
[\frac{g(I_0(t)e^{-k_0\rho(t)y})}{\mathcal{R}_0}\\[2mm]
\qquad  \; \; \; -d(\rho(t)y,t)-\frac{\dot{\rho}(t)}{\rho(t)}]w, \;& 0<y<L_0, \; 0<t\leq T,\\
w_y(y,t)=0,\;& y=0,L_0, \;  0<t\leq T,\\
w(y,0)=w(y,T),\;&  0<y<L_0,
\end{array}\right.
\label{b10}
\ees
it admits the same principal eigen-pair $(\mathcal{R}_0, w)$.

Now we multiply \eqref{b10} by $e^{\frac{\alpha}{D}\rho(t)y}\tilde{z}$, \eqref{b07} by $w$, and subtract the resulting equations to obtain
\bes
\begin{array}{llllll}
-[e^{\frac{\alpha}{D}\rho(t)y}\tilde{z} w_{t}+(e^{\frac{\alpha}{D}\rho(t)y}\tilde{z})_{t}w]&\leq& \frac{D}{\rho^2(t)}[(e^{\frac{\alpha}{D}\rho(t)y}w_y)_{y}\tilde{z}-(e^{\frac{\alpha}{D}\rho(t)y}\tilde{z}_y)_{y}w]\\[2mm]
&&+e^{\frac{\alpha}{D}\rho(t)y}[\frac{g(I_0(t)e^{-k_0\rho(t)y})}{\mathcal{R}_0}-\frac{g(I_0(t)e^{-k_0\rho(t)y})}{\mathcal{R}_1}]w\tilde{z}.
\end{array}
\label{b11}
\ees
Integrating the above inequality over $(0,L_0)\times(0,T)$ and using the initial boundary conditions of \eqref{b07}, \eqref{b10}, one has
\bes
\begin{array}{llllll}
&&\int_{0}^{L_0}\int_{0}^{T}e^{\frac{\alpha}{D}\rho(t)y}[\frac{g(I_0(t)e^{-k_0\rho(t)y})}{\mathcal{R}_0}-\frac{g(I_0(t)e^{-k_0\rho(t)y})}{\mathcal{R}_1}]w\tilde{z} dtdy\\[2mm] &\geq&-\int_{0}^{L_0}\int_{0}^{T}[e^{\frac{\alpha}{D}\rho(t)y}\tilde{z} w_{t}+(e^{\frac{\alpha}{D}\rho(t)y}\tilde{z})_{t}w]dtdy\\[2mm]
&&-\int_{0}^{L_0}\int_{0}^{T}\frac{D}{\rho^2(t)}[(e^{\frac{\alpha}{D}\rho(t)y}w_y)_{y}\tilde{z}-(e^{\frac{\alpha}{D}\rho(t)y}\tilde{z}_y)_{y}w]dtdy\\[2mm]
&=&-\int_{0}^{L_0}\int_{0}^{T}[e^{\frac{\alpha}{D}\rho(t)y}\tilde{z}w]_{t}dtdy-\int_{0}^{L_0}\int_{0}^{T}\frac{D}{\rho^2(t)}[(e^{\frac{\alpha}{D}\rho(t)y}w_y)\tilde{z}
-(e^{\frac{\alpha}{D}\rho(t)y}\tilde{z}_y)w]_{y}dtdy\\[2mm]
&=&0-\int_{0}^{T}\frac{D}{\rho^2(t)}[(e^{\frac{\alpha}{D}\rho(t)y}w_y)\tilde{z}-(e^{\frac{\alpha}{D}\rho(t)y}\tilde{z}_y)w]|_0^{L_0}dt\\[2mm]
&=&0,
\end{array}
\label{b12}
\ees
where the periodicity of $w$ and $z$ has been used. Moreover, recalling the hypothesis of $g$ and $I$ yields $\mathcal{R}_0\leq \mathcal{R}_1$. By \cite{ELC}, the positive eigenfunction corresponding to the principal eigenvalue $\mathcal{R}_0$ is unique ( subjects to a constant multiple ). Based on the above argument, we are confident that the last claim of lemma can be proved.
\epf

\bigskip
Next we consider the general case that $\alpha\neq0$. Multiplying both side of \eqref{b05} by $\frac{1}{\varphi}$, and integrating over $(0,L_0)\times(0,T)$ yield
\bes
\begin{array}{llllll}
0&=&\int_{0}^{T}\int_{0}^{L_0} [\frac{D}{\rho^2(t)} \frac{\varphi_{yy}}{\varphi}+\frac{\alpha}{\rho(t)}\frac{\varphi_y}{\varphi}+\frac{g(I_0(t)e^{-k_0\rho(t)y})}{\mathcal{R}_0}-d(\rho(t)y,t)]dydt\\[2mm]
&=&\int_{0}^{T}\int_{0}^{L_0} \frac{D}{\rho^2(t)} \frac{d\varphi_{y}}{\varphi}dt+\int_{0}^{T}\int_{0}^{L_0} \frac{\alpha}{\rho(t)}\frac{\varphi_y}{\varphi}dydt\\[2mm]
&&+\int_{0}^{T}\int_{0}^{L_0}[\frac{g(I_0(t)e^{-k_0\rho(t)y})}{\mathcal{R}_0}-d(\rho(t)y,t)]dydt\\[2mm]
&=&\int_{0}^{T}\int_{0}^{L_0} [\frac{D}{\rho^2(t)} \frac{\varphi_{y}^2}{\varphi^2}+\frac{\alpha}{\rho(t)}\frac{\varphi_y}{\varphi}]dydt\\[2mm]
&&+\int_{0}^{T}\int_{0}^{L_0}[\frac{g(I_0(t)e^{-k_0\rho(t)y})}{\mathcal{R}_0}
-d(\rho(t)y,t)]dydt\\[2mm]
&=&\int_{0}^{T}\int_{0}^{L_0} [D (\frac{\varphi_{y}}{\rho(t)\varphi}+\frac{\alpha}{2D})^2-\frac{\alpha^2}{4D}]dydt\\[2mm]
&&+\int_{0}^{T}\int_{0}^{L_0}[\frac{g(I_0(t)e^{-k_0\rho(t)y})}{\mathcal{R}_0}
-d(\rho(t)y,t)]dydt\\[2mm]
&\geq&\int_{0}^{T}\int_{0}^{L_0} (-\frac{\alpha^2}{4D})dydt+\int_{0}^{T}\int_{0}^{L_0}[\frac{g(I_0(t)e^{-k_0\rho(t)y})}{\mathcal{R}_0}
-d(\rho(t)y,t)]dydt,\\[2mm]
\end{array}
\label{b13}
\ees
which derive the following inequality
\bes
\begin{array}{llllll}
\mathcal{R}_0\geq\frac{\int_{0}^{T}\int_{0}^{L_0} g(I_0(t)e^{-k_0\rho(t)y})dydt}{\int_{0}^{T}\int_{0}^{L_0} [\frac{\alpha^2}{4D}+d(\rho(t)y,t)]dydt}.
\end{array}
\label{b14}
\ees
\begin{rmk}
It follows from \eqref{b13} that
\bes
\begin{array}{llllll}
\lim\limits_{D\rightarrow0,\alpha\rightarrow0}\mathcal{R}_0(D,\alpha, L)=\frac{\int_{0}^{T}\int_{0}^{L_0} g(I_0(t)e^{-k_0\rho(t)y})dydt}{\int_{0}^{T}\int_{0}^{L_0} d(\rho(t)y,t)dydt}:=\mathcal{R}_0^*.
\end{array}
\label{b31}
\ees
\end{rmk}

Now we consider the species case $\alpha=0$. Noting that
$g(I_0(t)e^{-k_0\rho(t)y})-d(\rho(t)y,t)$ is strictly decreasing with respect to $y$, which is straightforward to show that
\bes
\begin{array}{llllll}
&&\frac{D}{\rho^2(t)} \varphi_{yy}+[ \frac{g(I_0(t)e^{-k_0\rho(t)L_0})}{\mathcal{R}_0}-d(\rho(t)L_0,t)-\frac{\dot{\rho}(t)}{\rho(t)}]\varphi\\[2mm]
&\leq& \frac{D}{\rho^2(t)} \varphi_{yy}+[\frac{g(I_0(t)e^{-k_0\rho(t)y})}{\mathcal{R}_0}-d(\rho(t)y,t)-\frac{\dot{\rho}(t)}{\rho(t)}]\varphi \\[2mm]
&\leq&\frac{D}{\rho^2(t)} \varphi_{yy}+[\frac{g(I_0(t))}{\mathcal{R}_0}- d(0,t)-\frac{\dot{\rho}(t)}{\rho(t)}]\varphi.
\end{array}
\label{b15}
\ees

We then consider the following auxiliary problem
\bes\left\{
\begin{array}{llllll}
\hat{z}_{t}=\frac{D}{\rho^2(t)} \hat{z}_{yy}+
[\frac{g(I_0(t)e^{-k_0\rho(t)L_0})}{\mathcal{R}_2}-d(\rho(t)L_0,t)-\frac{\dot{\rho}(t)}{\rho(t)}]\hat{z}, \;& 0<y<L_0, \; 0<t\leq T,\\[2mm]
\hat{z}_y(y,t)=0,\;\;&y=0,L_0,\; 0<t\leq T, \\
\hat{z}(y,0)=\hat{z}(y,T),\;& 0<y<L_0.
\end{array}\right.
\label{b16}
\ees
Letting $\hat{z}(y,t)=f(t)$, \eqref{b16} turns to
\bes\left\{
\begin{array}{llllll}
f_{t}=[\frac{g(I_0(t)e^{-k_0\rho(t)L_0})}{\mathcal{R}_2}-d(\rho(t)L_0,t)-\frac{\dot{\rho}(t)}{\rho(t)}]f(t), \;&  0<t\leq T,\\
f(0)=f(T),
\end{array}\right.
\label{b17}
\ees
integrating \eqref{b17} over $(0,T)$ with a view to periodicity of $\rho(t)$ yields
\bes
\begin{array}{llllll}
\mathcal{R}_2=\frac{\int_{0}^{T} g(I_0(t))e^{-k_0\rho(t)L_0})dt}{\int_{0}^{T}d(\rho(t)L_0, t)dt}.
\end{array}
\label{b18}
\ees
It is accessible to obtain from \eqref{b16}
\bes\left\{
\begin{array}{llllll}
\hat{z}_{t}\leq\frac{D}{\rho^2(t)} \hat{z}_{yy}+
[\frac{g(I_0(t)e^{-k_0\rho(t)y})}{\mathcal{R}_2}-d(\rho(t)y,t)-\frac{\dot{\rho}(t)}{\rho(t)}]\hat{z}, \;& 0<y<L_0, \; 0<t\leq T,\\
\hat{z}_y(y,t)=0, \;&y=0,L_0,\; 0<t\leq T, \\
\hat{z}(y,0)=\hat{z}(y,T),\;& 0<y<L_0,
\end{array}\right.
\label{b19}
\ees
which together with Lemma 2.2 gives that $\mathcal{R}_2\leq\mathcal{R}_0$.

Similarly, we continue in this fashion to obtain
$\mathcal{R}_1\geq\mathcal{R}_0$
with
\bes
\begin{array}{llllll}
\mathcal{R}_1=\frac{\int_{0}^{T} g(I_0(t))dt}{\int_{0}^{T}d(0,t)dt}.
\end{array}
\label{b20}
\ees
We then conclude that
\bes
\begin{array}{ll}
\frac{\int_{0}^{T} g(I_0(t)e^{-k_0\rho(t)L_0})dt}{\int_{0}^{T} d(\rho(t)L_0,t)dt}
\leq\mathcal{R}_0(D,0,L)\leq\frac{\int_{0}^{T}g(I_0(t))dt}{\int_{0}^{T}d(0,t)dt}.
\end{array}
\label{b21}
\ees
More specifically, assume that most of the light is absorbed dependent of chlorophyll concentration \cite{SO,ky} (i.e. $k_0=0$) and $d(\rho(t)y,t)\equiv d(t)$, then $\mathcal{R}_0(D,0,L)=\frac{\int_{0}^{T}g(I_0(t))dt}{\int_{0}^{T}d(t)dt}$.

\bigskip
Now we are more interested in the relationship of $\mathcal{R}_0$ with respect to the the water column depth $L$ and the vertical turbulent diffusion rate $D$.
We first present the following monotonicity, whose proof is the same as that of \cite{pzh1} and so is omitted.
\begin{lem}
Assume that $\alpha=0$, for any given $D, L>0$, $\varphi_y(y,t)<0$ in $(0, L)$ for all $t$, where $\varphi(y,t)$ is the nonnegative eigenvalue function
corresponding to the principal eigenvalue of problem \eqref{b05}.
\end{lem}

Based on the lemma, we have the monotonicity of $\mathcal{R}_0$ with respect to the the water column depth $L$.
\begin{lem}
Assume that $\alpha=0$, for any given $D>0$, $\mathcal{R}_0(D,0,L)$ is strictly monotone decreasing in $L$.
\end{lem}
\bpf
For ease of notations, we write $\mathcal{R}_0(L)=\mathcal{R}_0(D,0,L)$. Given $0<L_0<L_1$, our goal is to show that $\mathcal{R}_0(L_0)> \mathcal{R}_0(L_1)$. We denote by $(\mathcal{R}_0(L_1), \varphi)$ and $(\mathcal{R}_0(L_0), w)$ are the principal eigen-pairs of \eqref{b06} and \eqref{b10}, respectively. By virtue of our notations, we have
\bes\left\{
\begin{array}{llllll}
\varphi_{t}=\frac{D}{\rho^2(t)}\varphi_{yy}+
[\frac{g(I_0(t)e^{-k_0\rho(t)y})}{\mathcal{R}_0(L_1)}
-d(\rho(t)y,t)-\frac{\dot{\rho}(t)}{\rho(t)}]\varphi, \;& 0<y<L_1, \; 0<t\leq T,\\
\varphi_y(y,t)=0,\;& y=0,L_1, \;  0<t\leq T,\\
\varphi(y,0)=\varphi(y,T),\;&  0<y<L_1,
\end{array}\right.
\label{b22}
\ees
and
\bes\left\{
\begin{array}{llllll}
-w_{t}=\frac{D}{\rho^2(t)}w_{yy}+
[\frac{g(I_0(t)e^{-k_0\rho(t)y})}{\mathcal{R}_0(L_0)}
-d(\rho(t)y,t)-\frac{\dot{\rho}(t)}{\rho(t)}]w, \;& 0<y<L_0, \; 0<t\leq T,\\
w_y(y,t)=0,\;& y=0,L_0, \;  0<t\leq T,\\
w(y,0)=w(y,T),\;&  0<y<L_0.
\end{array}\right.
\label{b23}
\ees
Now multiplying \eqref{b22} by $w$, and \eqref{b23} by $\varphi$, then subtracting the resulting equations to obtain
\bes
\begin{array}{llllll}
-[\varphi w_{t}+\varphi_{t}w]&=& \frac{D}{\rho^2(t)}[w_{yy}\varphi-\varphi_{yy}w]
+\varphi w[\frac{g(I_0(t)e^{-k_0\rho(t)y})}{\mathcal{R}_0(L_0)}-\frac{g(I_0(t)e^{-k_0\rho(t)y})}
{\mathcal{R}_0(L_1)}].
\end{array}
\label{b24}
\ees
Integrating the above equation over $(0,L_0)\times(0,T)$ and using the initial boundary conditions of \eqref{b22}, \eqref{b23}, one has
\bes
\begin{array}{llllll}
&&\int_{0}^{L_0}\int_{0}^{T}\varphi w[\frac{g(I_0(t)e^{-k_0\rho(t)y})}{\mathcal{R}_0(L_0)}-\frac{g(I_0(t)e^{-k_0\rho(t)y})}{\mathcal{R}_0(L_1)}] dtdy\\[2mm]
&=&-\int_{0}^{L_0}\int_{0}^{T}[\varphi w_{t}+\varphi_{t}w]dtdy-\int_{0}^{T}\int_{0}^{L_0}\frac{D}{\rho^2(t)}[w_y\varphi-\varphi_y w]_{y}dy dt\\[2mm]
&=&\int_{0}^{T}\frac{D}{\rho^2(t)}w_y\varphi_y(L_0,t)w(L_0,t)dt\\[2mm]
&<&0,
\end{array}
\label{b25}
\ees
in which we have used the fact that $\varphi_y(y,t)<0$ for $0<y<L_1, 0<t<T$ due to Lemma 2.3, we then have $\mathcal{R}_0(L_0)>\mathcal{R}_0(L_1)$.
\epf
\begin{lem}
Assume that $\alpha=0$, for any given $D>0$, $\mathcal{R}_0(D,0,L)$ is strictly monotone decreasing in $D$.
\end{lem}
\bpf
Let $\mathcal{R}_0=\mathcal{R}_0(D,0,L)$ and $(\mathcal{R}_0, \varphi)$ is the principal eigen-pair of \eqref{b05}. Clearly, $\mathcal{R}_0(D,0,L)$ and $\varphi$ are $C^1-$ function of $D$. For simplicity of presentation, we denote $\frac{\partial\varphi}{\partial D}$ by $\varphi^{'}$ and $\frac{\partial \mathcal{R}_0}{\partial D}$ by $\mathcal{R}_0^{'}$. Letting $\alpha=0$, by assumption, and differentiating \eqref{b05} with respect to $D$ yields
{\small \bes\left\{
\begin{array}{lll}
\varphi^{'}_{t}=\frac{D}{\rho^2(t)} \varphi^{'}_{yy}+\frac{1}{\rho^2(t)}\varphi_{yy}+[\frac{g(I_0(t)e^{-k_0\rho(t)y})}{\mathcal{R}_0}-d(\rho(t)y,t)-\frac{\dot{\rho}(t)}{\rho(t)}]\varphi^{'}\\[2mm]
\qquad-(\frac{1}{\mathcal{R}_0})^{'}g(I_0(t)e^{-k_0\rho(t)y})\varphi,& 0<y<L_0, \; 0<t\leq T,\\[2mm]
\varphi^{'}_y(y,t)=0,\;& y=0,L_0, \;  0<t\leq T,\\[2mm]
\varphi^{'}(y,0)=\varphi^{'}(y,T),\;& 0<y<L_0,
\end{array}\right.
\label{b26}
\ees}
Let $w$ be the principal eigenfunction corresponding to $\mathcal{R}_0$ which satisfies \eqref{b10}. Multiplying \eqref{b26} by $w$ and integrating the equality over $(0,L_0)\times(0,T)$, by integration by parts, yields
{\small
\bes\left\{
\begin{array}{lll}
-\int_0^{L_0}\int_0^{T}\varphi^{'}w_{t}dtdy=-D\int_0^{L_0}\int_0^{T}\frac{1}{\rho^2(t)}\varphi^{'}_{y}w_y dtdy\\[2mm]
\quad +\int_0^{L_0}\int_0^{T}w\Big\{\frac{1}{\rho^2(t)}\varphi_{yy}
+[\frac{g(I_0(t)e^{-k_0\rho(t)y})}{\mathcal{R}_0}-d(\rho(t)y,t)-\frac{\dot{\rho}(t)}{\rho(t)}]\varphi^{'}\\[2mm]
\quad -(\frac{1}{\mathcal{R}_0})^{'}g(I_0(t)e^{-k_0\rho(t)y})\varphi\Big\}dtdy, &0<y<L_0, \; 0<t\leq T,\\[2mm]
\varphi^{'}_y(y,t)=0,& y=0, L_0,  \; 0<t\leq T,\\[2mm]
\varphi^{'}(y,0)=\varphi^{'}(y,T),&0<y<L_0.
\end{array}\right.
\label{b27}
\ees
}
By substituting $-w_{t}=\frac{D}{\rho^2(t)}w_{yy}+
[\frac{g(I_0(t)e^{-k_0\rho(t)y})}{\mathcal{R}_0}
-d(\rho(t)y,t)-\frac{\dot{\rho}(t)}{\rho(t)}]w$ into \eqref{b27}, we get
\bes
\begin{array}{lll}
(\frac{1}{\mathcal{R}_0})^{'}=-\frac{\int_0^{L_0}\int_0^{T}\frac{1}{\rho^2(t)}\varphi_{yy}wdtdy}
{\int_0^{L_0}\int_0^{T}g(I_0(t)e^{-k_0\rho(t)y})\varphi w dtdy}
=\frac{\int_0^{L_0}\int_0^{T}\frac{e^{\frac{\alpha}{D}\rho(t)y}}{\rho^2(t)}\varphi_y w_ydtdy}
{\int_0^{L_0}\int_0^{T}g(I_0(t)e^{-k_0\rho(t)y})\varphi w dtdy},
\end{array}
\label{b28}
\ees
Considering \eqref{b10}, set $\xi(y,t)=w(y,-t)$, we have
\bes\left\{
\begin{array}{lll}
\xi_{t}=\frac{D}{\rho^2(t)}\xi_{yy}+
[\frac{g(I_0(-t)e^{-k_0\rho(-t)y})}{\mathcal{R}_0}
-d(\rho(-t)y,t)-\frac{\dot{\rho}(-t)}{\rho(-t)}]\xi, \;& 0<y<L_0, \; 0<t\leq T,\\
\xi_y(y,t)=0,\;& y=0,L_0, \;  0<t\leq T,\\
\xi(y,0)=\xi(y,T),\;&  0<y<L_0,
\end{array}\right.
\label{b29}
\ees
In accordance to $d(\rho(t)y,t)$ is strictly decreasing with respect to $y$, using the same argument as in the proof of Lemma 2.3, we can easily carry out
\bes
\begin{array}{lll}
\xi_y(y,t)=w_y(y,-t)<0, \qquad  0<y<L_0, \; 0<t\leq T.
\end{array}
\label{b30}
\ees
So \eqref{b28} implies $(\frac{1}{\mathcal{R}_0})^{'}>0$.
\epf

\section{The threshold-type dynamics}

The following result shows the long time behavior of solution depends on the basic reproduction number $\mathcal{R}_0$.
\begin{thm}
The following statements are valid:

$(i)$ If  $\mathcal{R}_0\leq1$, then for any solution $v(y,t)$ of system  \eqref{a07}, \eqref{a08}, we can deduce that
$\lim\limits_{t\to\infty}{v(y,t)}=0$ for $y\in[0,L_0]$, that is, the phytoplankton species dies out in a long run.

$(ii)$ If $\mathcal{R}_0>1$, then problem
\eqref{a07}, \eqref{a09} admits at least one positive periodic solution. Moreover, there exists $\delta_0>0$ such that any positive solution $v(y,t)$ of system
\eqref{a07}, \eqref{a08} satisfies  $\limsup\limits_{t\to\infty}v(y,t)\geq\delta_0$ uniformly for $y\in [0,L_0]$, that is, the phytoplankton species persists uniformly.
\label{main}
\end{thm}
\bpf
$(i)$ When $\mathcal{R}_0<1$, problem \eqref{b03} admits an eigen-pair $(\mathcal{R}_0;\phi)$ such that $\phi(y,t)>0$ for $(y,t)\in[0,L_0]\times[0,T]$. Letting $\tilde{v}(y,t)=Me^{-\lambda t}\phi(y,t),$ where $0<\lambda\leq g(I_0e^{-k_0\rho(t)y})(\frac{1}{\mathcal{R}_0}-1)$ for $(y,t)\in[0,L_0]\times[0,T]$, from the hypothesis of $g$ and mathematical expression of $I$, we then have
$$\begin{array}{llllll}
&&\tilde{v}_{t}-\frac{D}{\rho^2(t)} \tilde{v}_{yy}+\frac{\alpha}{\rho(t)}\tilde{v}_y-[g(I(y,t))-d(y,t)-\frac{\dot{\rho}(t)}{\rho(t)}]\tilde{v}\\[2mm]
&\geq&\tilde{v}_{t}-\frac{D}{\rho^2(t)} \tilde{v}_{yy}+\frac{\alpha}{\rho(t)}\tilde{v}_y-[g(I_0e^{-k_0\rho(t)y})-d(y,t)-\frac{\dot{\rho}(t)}{\rho(t)}]\tilde{v}\\[2mm]
&=&Me^{-\lambda t}\phi_t-\lambda Me^{-\lambda t}\phi-Me^{-\lambda t}\frac{D}{\rho^2(t)} \phi_{yy}+Me^{-\lambda t}\frac{\alpha}{\rho(t)}\phi_y\\[2mm]
&&-Me^{-\lambda t}\phi[g(I_0e^{-k_0\rho(t)y})-d(y,t)]+ Me^{-\lambda t}\frac{\dot \rho(t)}{\rho(t)}\phi\\[2mm]
&=&\tilde{v}\{-\lambda+\frac{g(I_0e^{-k_0\rho(t)y})}{\mathcal{R}_0}-g(I_0e^{-k_0\rho(t)y})\}\\[2mm]
&\geq& 0,
\end{array}$$
for $(y,t)\in (0,L_0)\times(0,\infty)$ and
$$D\tilde{v}_y-\alpha\rho(t)\tilde{v}(y,t)=0,\; y=0,L_0, \;  t>0,$$
$$\tilde{v}(y,0)=u_0(y)\geq0,\not\equiv0, \;0<y<L_0,$$
therefore $\tilde{v}$ is the upper solution of \eqref{a07}, \eqref{a08}
if $M$ is large enough. Since $\lim\limits_{t\rightarrow+\infty}\tilde{v}(y,t)=0,$ then $\lim\limits_{t\rightarrow+\infty}v(y,t)=0$ uniformly for $y\in[0,L_0]$.

Next, we claim that \eqref{a07}, \eqref{a09} has no positive solution if $\mathcal{R}_0\leq1$. Assume that $v^*(y,t)$ is a positive solution of \eqref{a07}, \eqref{a09}. By the assumption that $g$ is strictly increasing with respect to $I$ and $I=I_0e^{-k_0\rho(t)y-k_1\int_0^{\rho(t)y}v^*(s,t)ds}$, we obtained
\bes\left\{
\begin{array}{lll}
v^*_{t}-\frac{D}{\rho^2(t)} v^*_{yy}+\frac{\alpha}{\rho(t)}v^*_y-[g(I_0(t)e^{-k_0\rho(t)y})\\[2mm]
\qquad -d(\rho(t)y,t)-\frac{\dot{\rho}(t)}{\rho(t)}]v^*\leq,\not\equiv 0, \;& 0<y<L_0, \; 0<t\leq T,\\
Dv^*_y-\alpha\rho(t)v^*(y,t)=0,\;& y=0,L_0, \;  0<t\leq T,\\
v^*(y,0)=v^*(y,T),\;& 0<y<L_0.
\end{array}\right.
\label{c01}
\ees
In a similar way, let $z(y,t)=e^{-\frac{\alpha}{D}\rho(t)y}v^*(y,t)$, thus $z>0$ and $z$ satisfies
\bes\left\{
\begin{array}{lll}
z_{t}-\frac{D}{\rho^2(t)} z_{yy}-\frac{\alpha}{\rho(t)}z_y-[g(I_0(t)e^{-k_0\rho(t)y})\\[2mm]
\qquad -d(\rho(t)y,t)-\frac{\dot{\rho}(t)}{\rho(t)}]z\leq,\not\equiv0 \;& 0<y<L_0, \; 0<t\leq T,\\
z_y(y,t)=0,\;& y=0,L_0, \;  0<t\leq T,\\
z(y,0)=z(y,T),\;& 0<y<L_0.
\end{array}\right.
\label{c02}
\ees
In connection with Proposition 5.2(ii) \cite{pzh}, which deserve $\lambda_0<0$. In addition, we deduce from Lemma 2.1(ii) that $\mathcal{R}_0>1$, which contradicts $\mathcal{R}_0\leq1$. So we admit that \eqref{a07}, \eqref{a09} does not have any positive solution provided that $\mathcal{R}_0\leq1$.

For $\mathcal{R}_0=1$, then $\lambda_0=0$, \eqref{a07}, \eqref{a09} has no positive solution, it is pretty easy for us to verify conclusion (i) by making use of Theorem 3.4 \cite{DH} and Lemma
2.4 \cite{pzh}.

$(ii)$ When $\mathcal{R}_0>1$,
assume, for the sake of contradiction, that for any $\delta$ there exists a positive solution $v(y,t)$ of problem \eqref{a07} and \eqref{a08} such that
\begin{equation}
\limsup\limits_{t\to\infty}v(y,t)\leq\frac{\delta}{2}.
\label{c03}
\end{equation}
For the above given $\delta$, there exists $T_{\delta}$ such that
$$0\leq v(y,t)\leq\delta \ \textrm{for} \;(y,t)\in [0,L_0]\times [T_{\delta}, +\infty).$$
Since $I(\rho(t)y,t)=I_0e^{-k_0\rho(t)y-k_1\int_{0}^{\rho(t)y}v(s,t)ds}$, then
$$I_0e^{-(k_0+k_1\delta)\rho(t)y}\leq I(\rho(t)y,t)\leq I_0e^{-k_0\rho(t)y}.$$
Then we have
\begin{equation}
\begin{array}{llllll}
v_{t}-\frac{D}{\rho^2(t)} v_{yy}+\frac{\alpha}{\rho(t)}v_y
&=&[g(I(\rho(t)y,t))-d(\rho(t)y,t)-\frac{\dot{\rho}(t)}{\rho(t)}]v\\
&\geq&[g(I_0e^{-(k_0+k_1\delta)\rho(t)y})-d(\rho(t)y,t)-\frac{\dot{\rho}(t)}{\rho(t)}]v
\end{array}
\label{c04}
\end{equation}
for $(y,t)\in [0,L_0]\times [T_{\delta}, +\infty)$. We now choose a sufficiently small number $\eta>0$
such that
\begin{equation}
v(y, T_{\delta})\geq \eta\phi(y, T_{\delta}),
\label{c05}
\end{equation}
where $\phi(y,t)>0$ for $(y,t)\in[0,L_0]\times [T_{\delta}, +\infty)$ satisfies \eqref{b03} with $\mathcal{R}_0>1$.
Set $0<\lambda_0\leq(1-\frac{1}{\mathcal{R}_0})g(I_0e^{-k_0\rho(t)y})$ for $y\in(0.L_0)$, and direct calculations show that, $\underline{v}(y,t)=\eta e^{\lambda_0(t-T_\delta)}\phi(y,t)$ satisfies
\begin{eqnarray}
\left\{
\begin{array}{lllll}
\underline v_{t}\leq\frac{D}{\rho^2(t)} \underline v_{yy}-\frac{\alpha}{\rho(t)}\underline v_y+[g(I_0e^{-(k_0+k_1\delta)\rho(t)y})\\[2mm]
\qquad\ -d(\rho(t)y,t)
-\frac{\dot{\rho}(t)}{\rho(t)}]\underline v, \;& 0<y<L_0, \; 0<t\leq T_\delta,\\
D\underline v_y-\alpha\rho(t)\underline v(y,t)=0,\ \;& y=0,L_0, \;  0<t\leq T_\delta,\\[2mm]
\underline v(y, T_{\delta}) =\eta\phi(y,T_\delta)\leq v(y, T_{\delta}), \;& 0<y<L_0.
\end{array} \right.
\label{c06}
\end{eqnarray}
It follows from \eqref{c04} and the comparison principle that
$$v(y,t)\geq \underline v(y, t)=\eta e^{\lambda_0(t-T_\delta)}\phi(y,t) \;\rm{for}\;   \textit{y}\in[0,L_0],\ \textit{t}\geq T_{\delta},$$
therefore, $v(y,t) \rightarrow\infty$ as $t\rightarrow\infty$, which contradicts \eqref{c03}.
This proves statement $(ii)$.
\epf
\section{\bf Simulation and discussion}

\begin{figure}[ht]
\centering
\subfigure[]{ {
\includegraphics[width=0.3\textwidth]{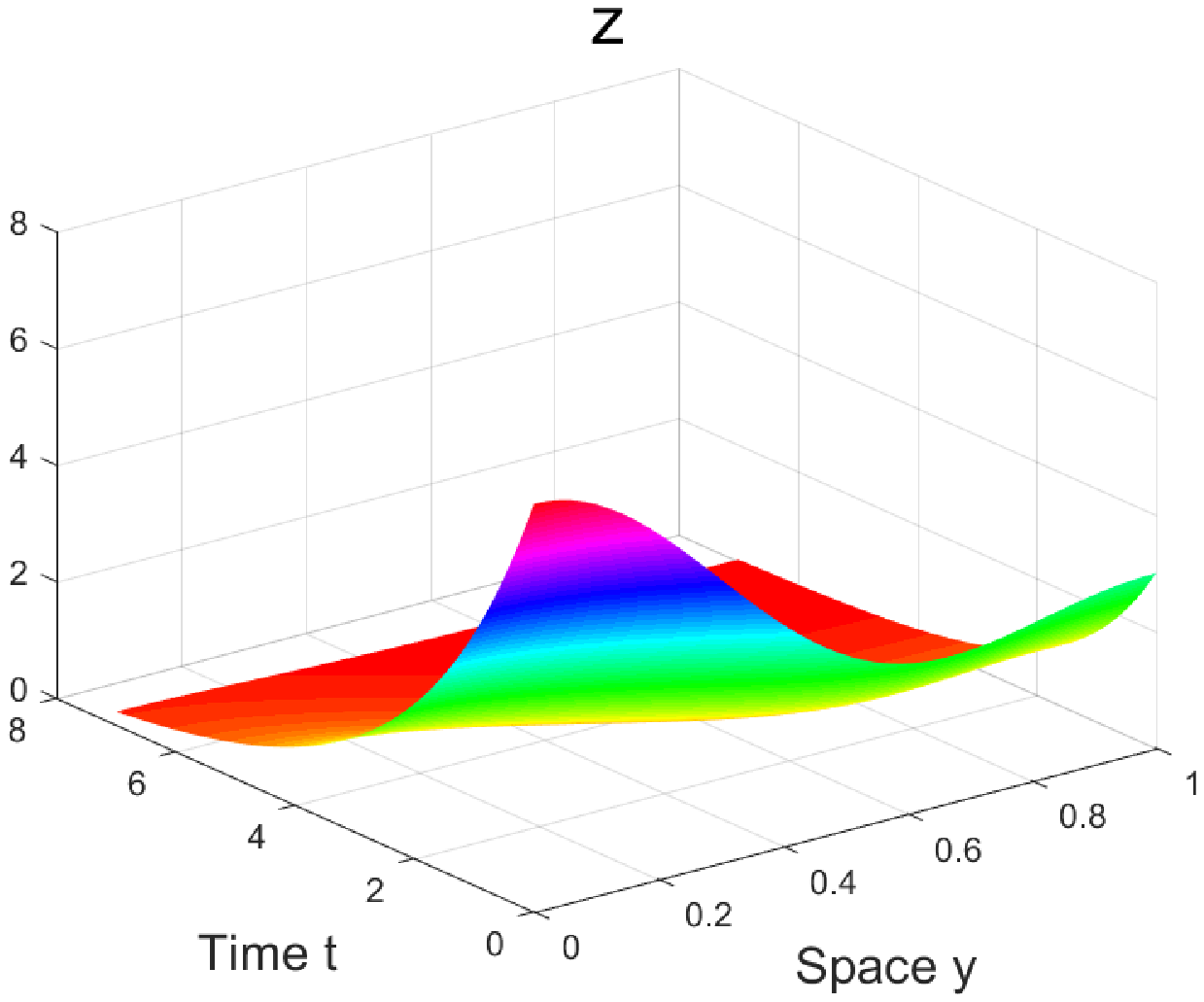}
} }
\subfigure[]{ {
\includegraphics[width=0.3\textwidth]{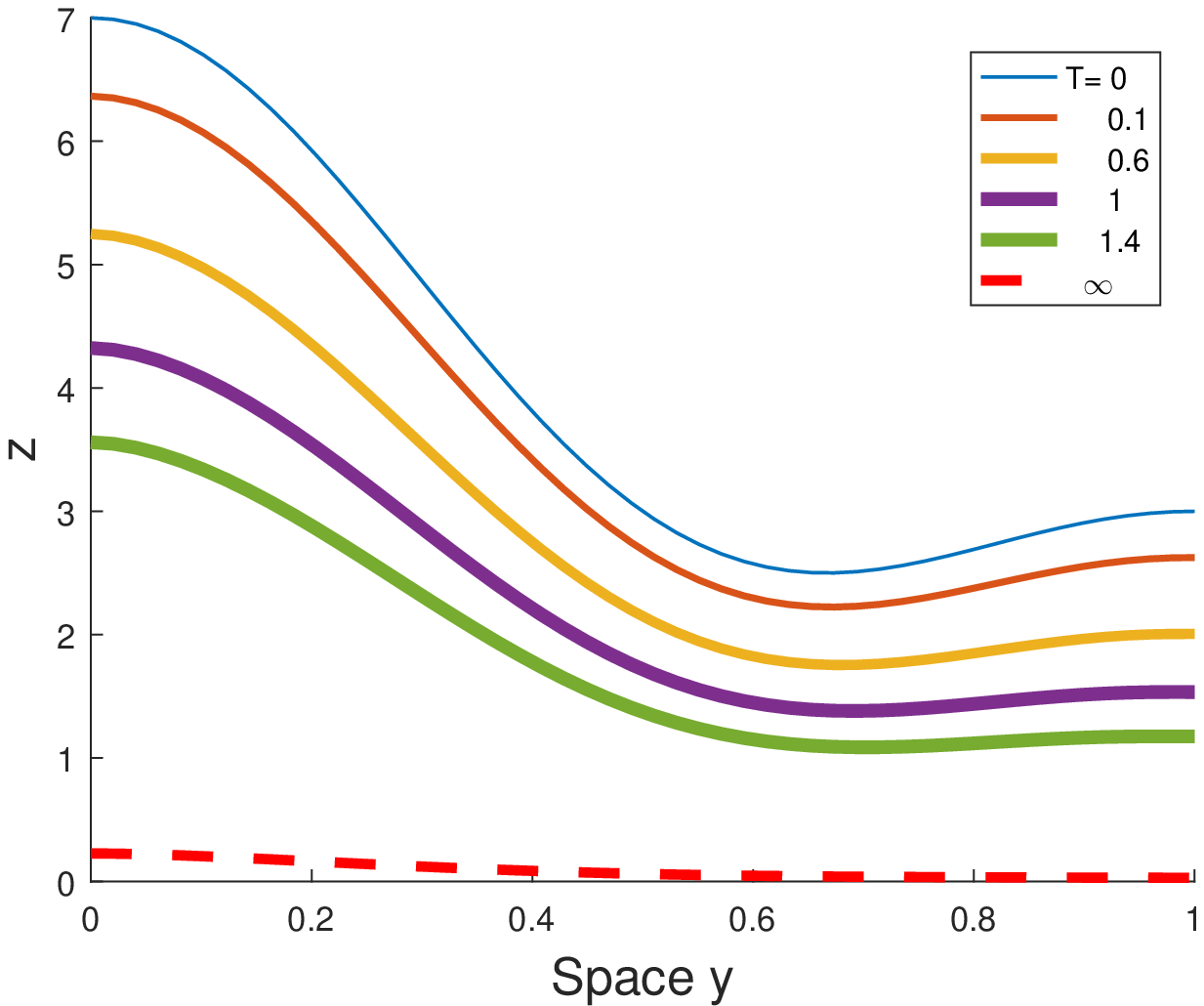}
} }
\subfigure[]{ {
\includegraphics[width=0.3\textwidth]{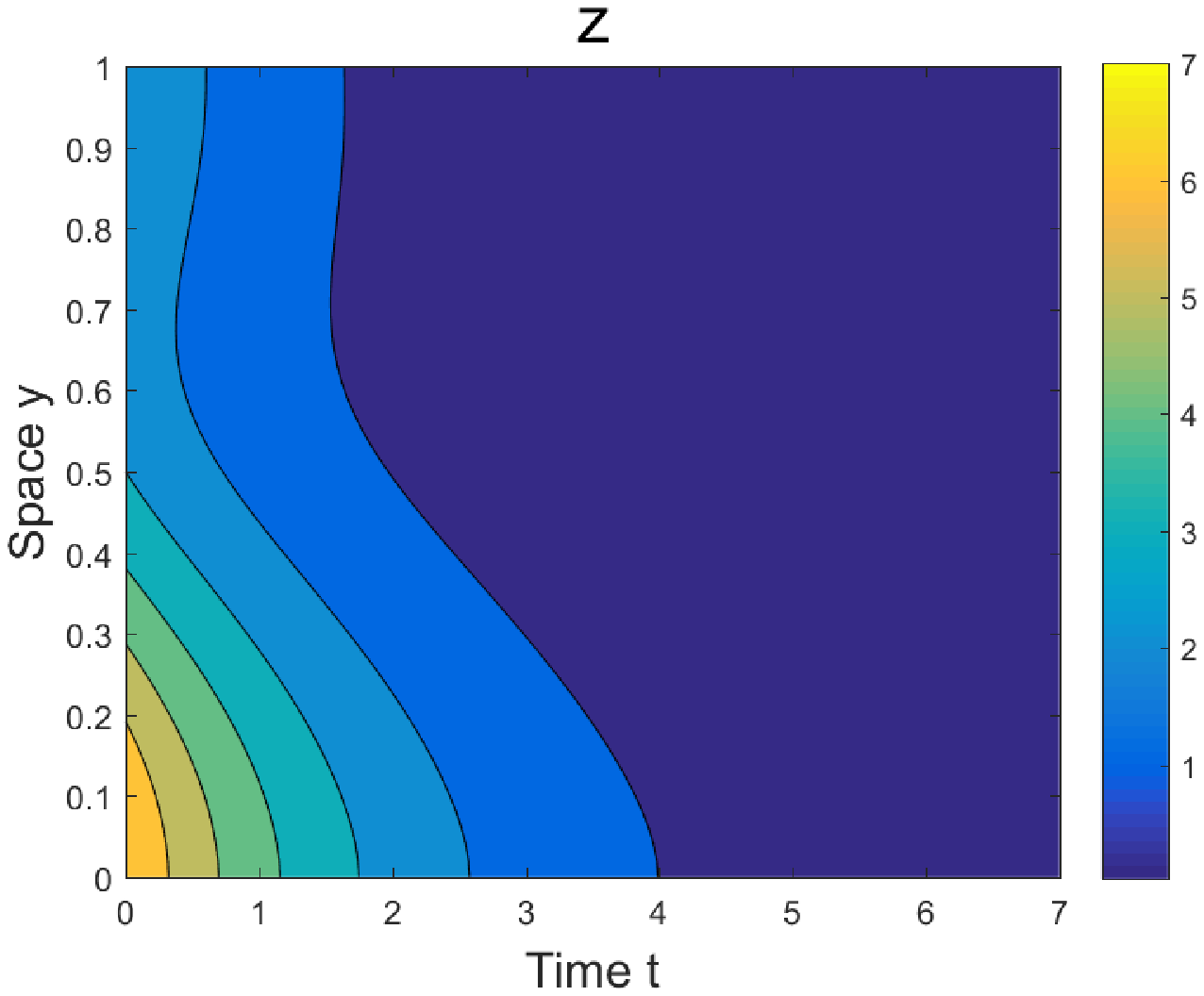}
} }
\caption{\scriptsize  $\rho_1(t)\equiv 1$. The interval is fixed and ${\mathcal{R}}_0<1$. Graphs $(a)-(c)$ showed that $z$ decay to $0$.
Graphs $(b)$ and $(c)$ are the cross-sectional view and contour map, respectively.}
\label{tu1}
\end{figure}

\begin{figure}[ht]
\centering
\subfigure[]{ {
\includegraphics[width=0.30\textwidth]{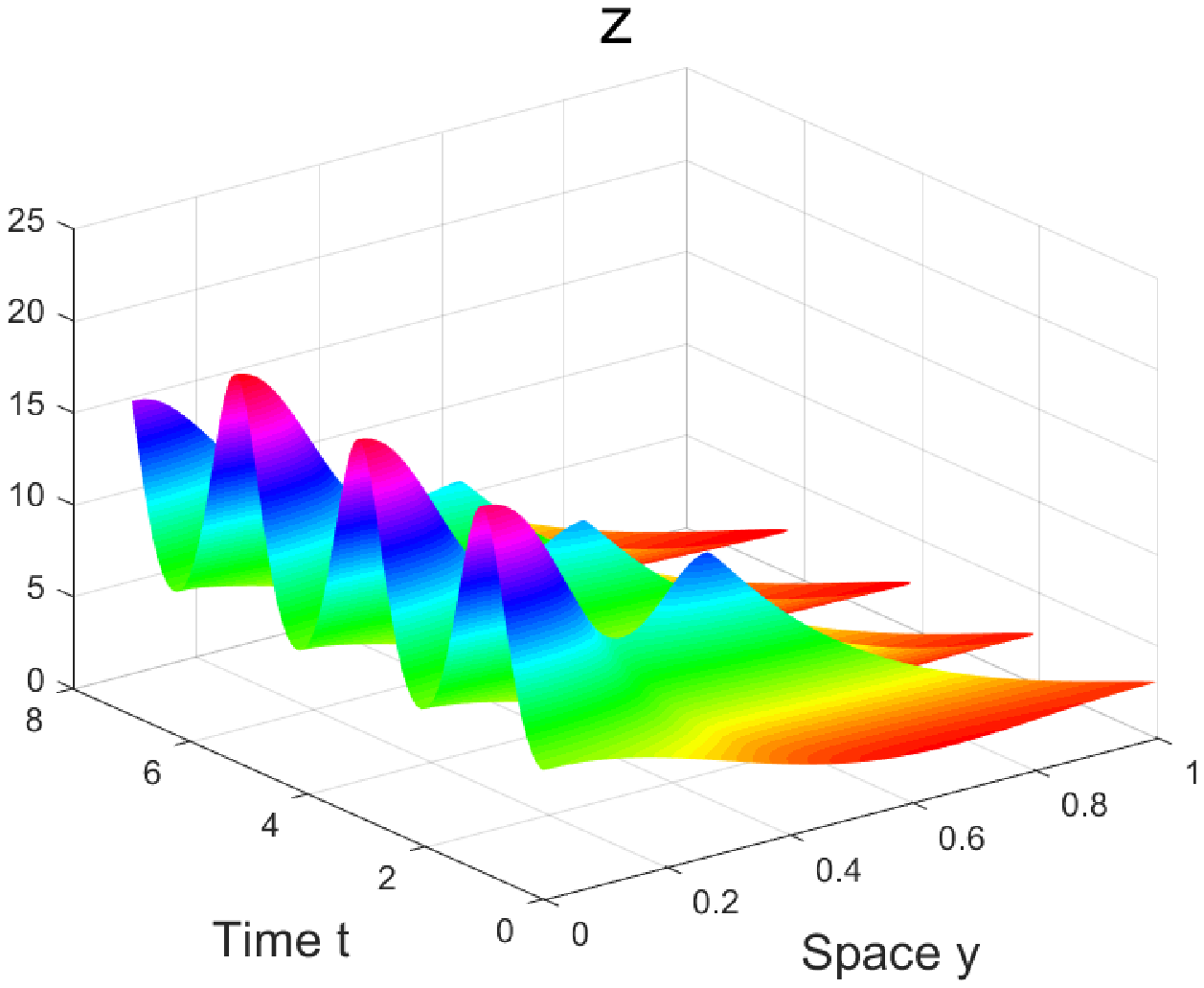}
} }
\subfigure[]{ {
\includegraphics[width=0.30\textwidth]{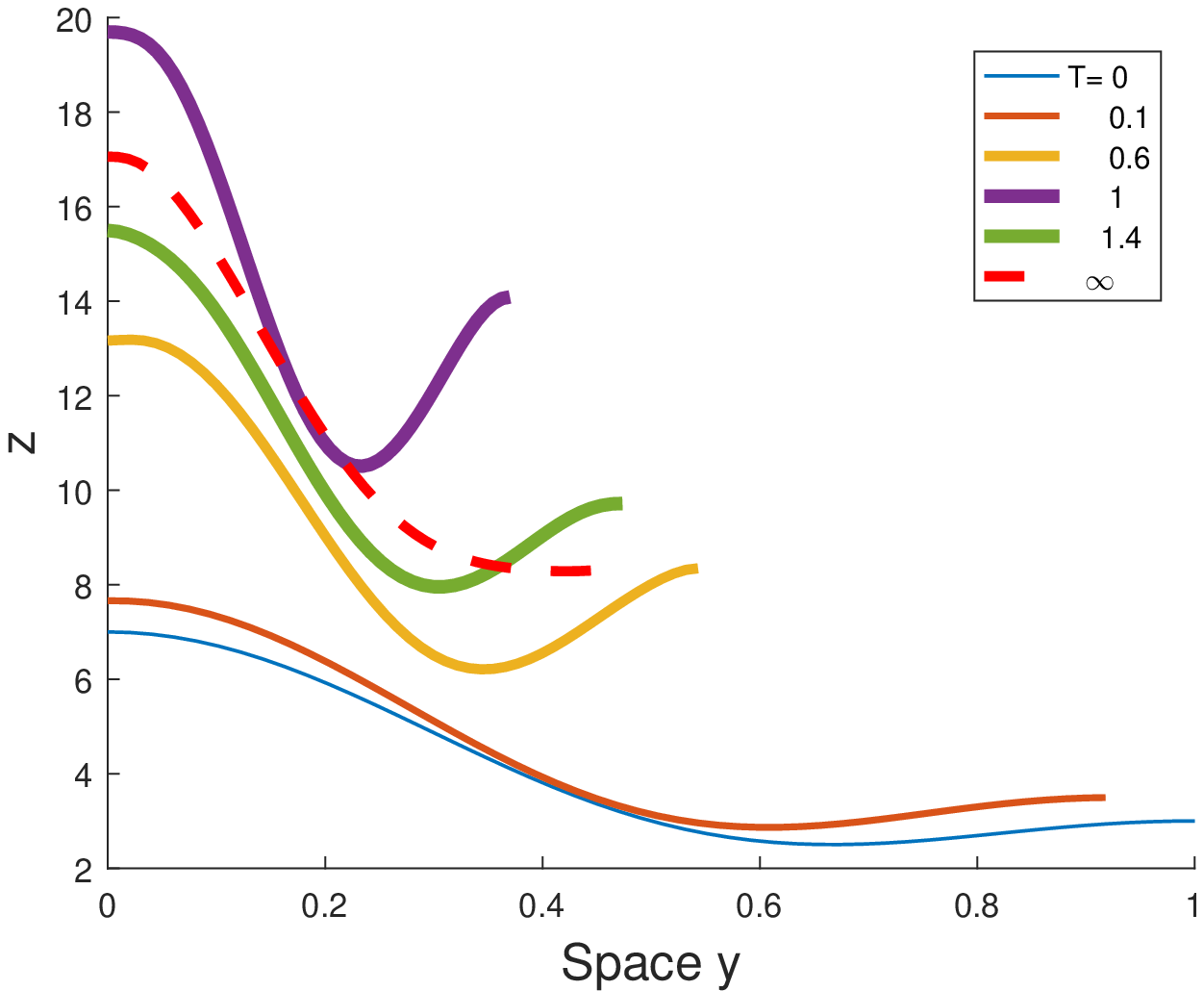}
} }
\subfigure[]{ {
\includegraphics[width=0.30\textwidth]{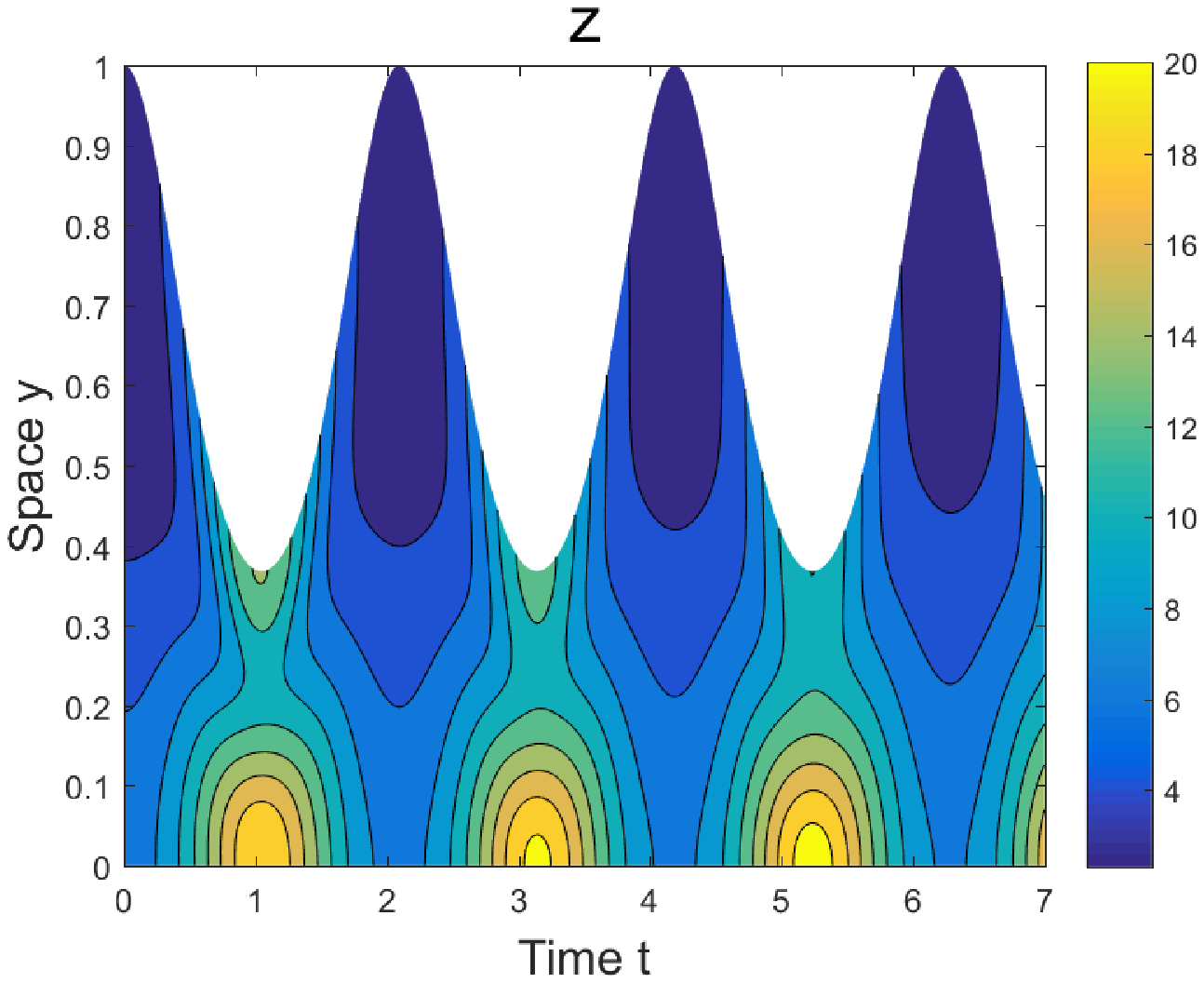}
} }
\caption{\scriptsize  $\rho_2(t)=e^{-0.5(1-\cos(3t))}$. The interval is evolving with a larger evolution rate $\rho_2(t)$ and ${\mathcal{R}}_0>1$. Graphs $(a)$ show that $z$ stabilized to a positive periodic steady state. Graphs $(b)$ and $(c)$, which are the cross-sectional view and contour map respectively, present the periodic
evolution of the interval.}
\label{tu2}
\end{figure}

In this section, we first carry out numerical simulations for problem \eqref{a07}, \eqref{a08} to illustrate the theoretical results by using Matlab.
The term $g(I(\rho(t)y,t))$ is the specific growth rate of phytoplankton as a function of light intensity $I(\rho(t)y,t)$. A standard model that incorporates saturation is $g(I)=\frac{aI}{1+bI}$ (see \cite{ky}).
Let us fix some coefficients. Assume that
\begin{equation*}
\left.
\begin{array}{lll}
D=0.001,\ \alpha=0.001, \  \ a=3,\ b=2,\   I=I_0e^{-k_0\rho(t)y}, \ I_0=0.1,\ k_0=0.2, \\[7pt]
L_0=1, \ z_0(y)=4+2\cos(\pi y)+\cos(2\pi y),
\end{array}
\right.
\end{equation*}
then the asymptotic behaviors of the solution to problem  are shown by choosing different $\rho(t)$ and $d(\rho(t)y,t)$.
\begin{exm}
Fix $d(\rho(t)y,t)=0.2+0.1\rho(t)y$. We first choose $\rho_1(t)\equiv 1$, which means that the habitat is fixed. From \eqref{b31}, direct calculations show that
$${\mathcal{R}}_0^*(\rho_1)=\frac{\int_{0}^{T}\int_{0}^{L_0} g(I_0(t)e^{-k_0\rho(t)y})dydt}{\int_{0}^{T}\int_{0}^{L_0} d(\rho(t)y,t)dydt}\approx\frac{0.4819}{0.5236}<1.$$
It is easily seen from Fig. \ref{tu1} that the phytoplankton species $z$ decays to zero.

Now we choose $\rho_2(t)=e^{-0.5(1-\cos(3t))}$, it follows from \eqref{b31} that
$${\mathcal{R}}_0^*(\rho_2)=\frac{\int_{0}^{T}\int_{0}^{L_0} g(I_0(t)e^{-k_0\rho(t)y})dydt}{\int_{0}^{T}\int_{0}^{L_0} d(\rho(t)y,t)dydt}
\approx\frac{0.4963}{0.4864}>1.$$
Easily to see from Fig. \ref{tu2} that the phytoplankton species $z$ stabilizes to a positive periodic steady state.

 One can see from the example that the phytoplankton vanishes in a fixed interval, but persist in a periodically evolving interval.
\end{exm}

\begin{figure}[ht]
\centering
\subfigure[]{ {
\includegraphics[width=0.30\textwidth]{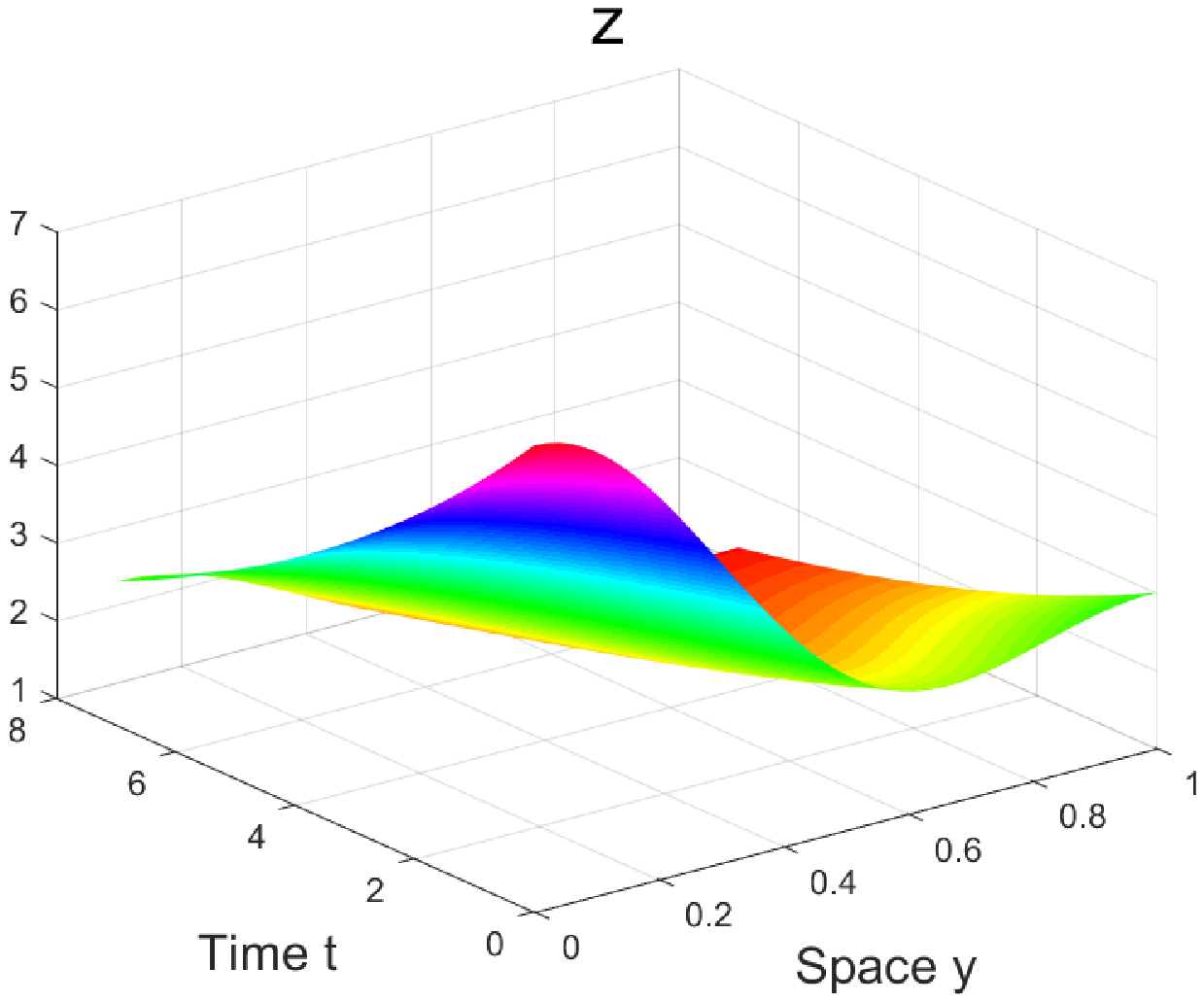}
} }
\subfigure[]{ {
\includegraphics[width=0.30\textwidth]{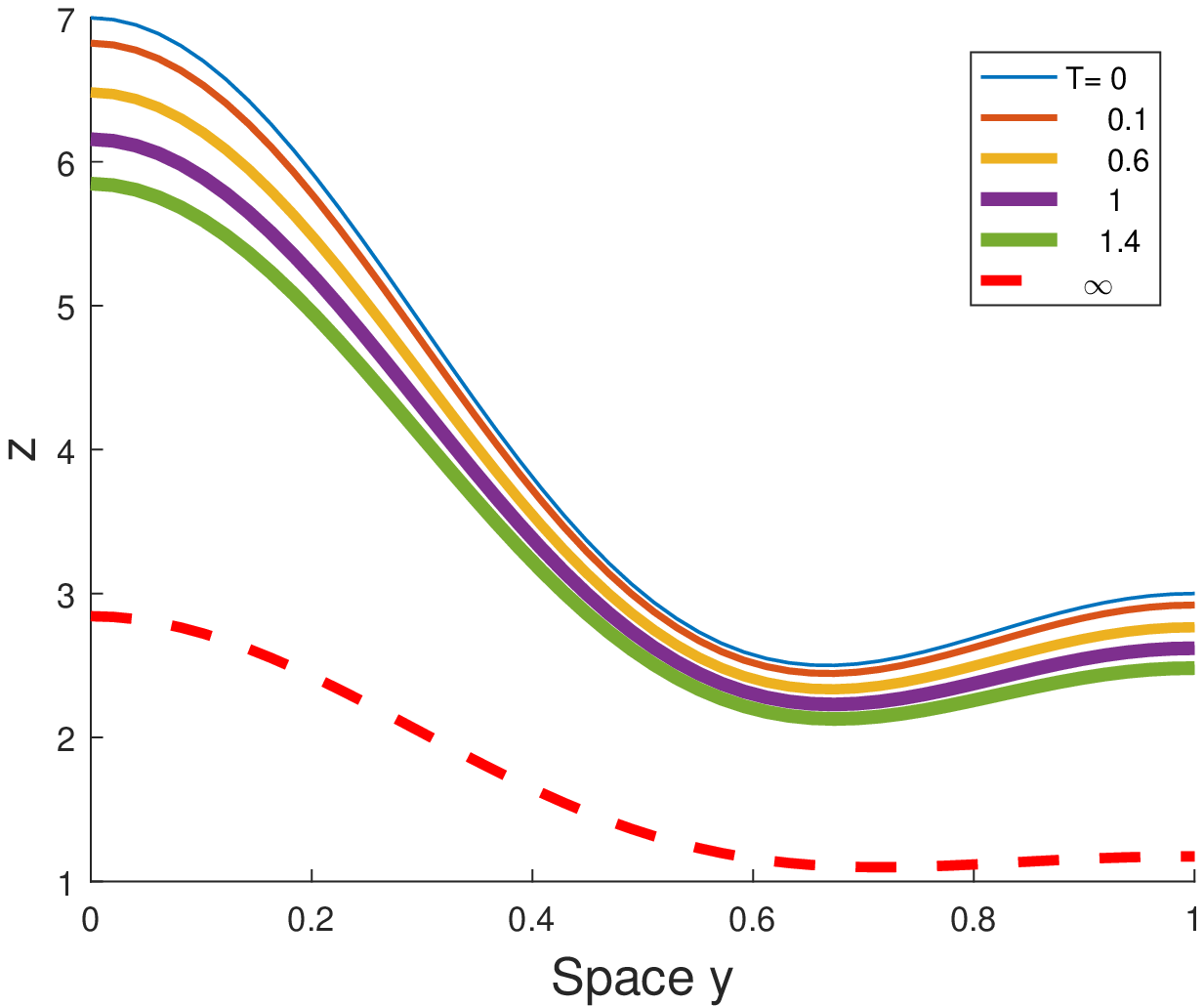}
} }
\subfigure[]{ {
\includegraphics[width=0.30\textwidth]{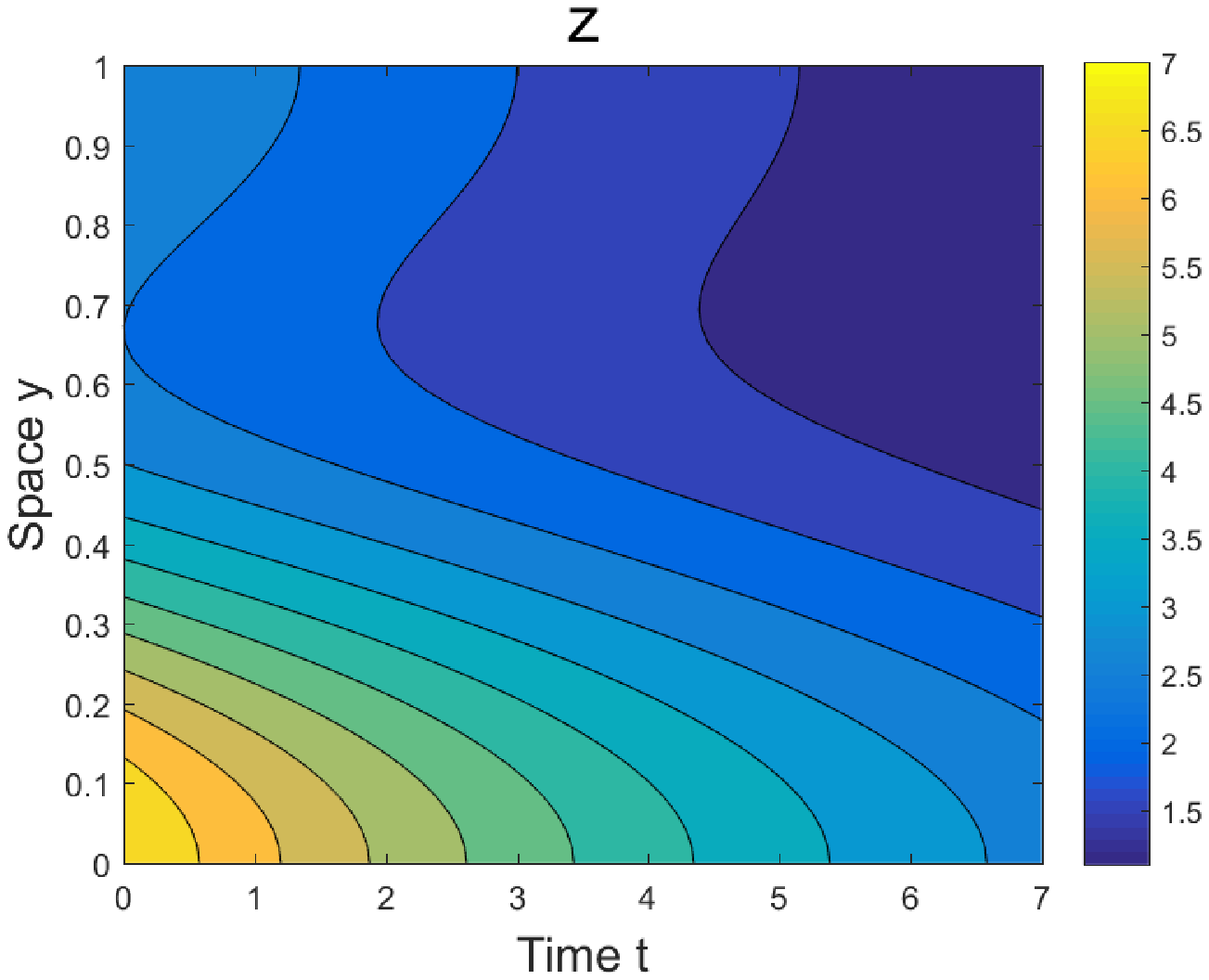}
} }
\caption{\scriptsize $\rho_3(t)\equiv 1$. The interval is fixed and ${\mathcal{R}}_0>1$. Graphs $(a)-(c)$ showed that $z$ decays to $0$. Graphs $(b)$ and $(c)$ are the cross-sectional view and contour map respectively.}
\label{tu3}
\end{figure}

\begin{figure}[ht]
\centering
\subfigure[]{ {
\includegraphics[width=0.30\textwidth]{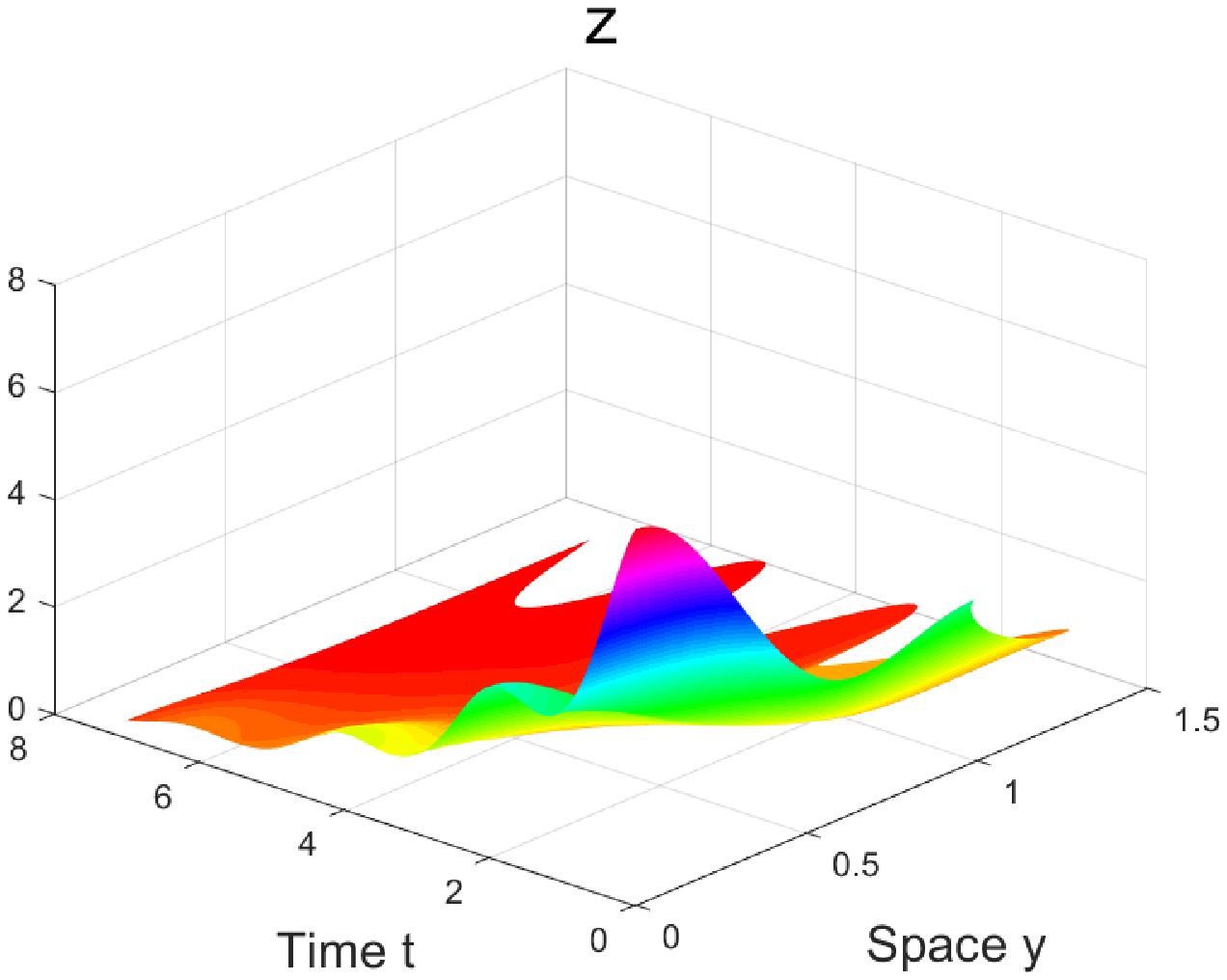}
} }
\subfigure[]{ {
\includegraphics[width=0.30\textwidth]{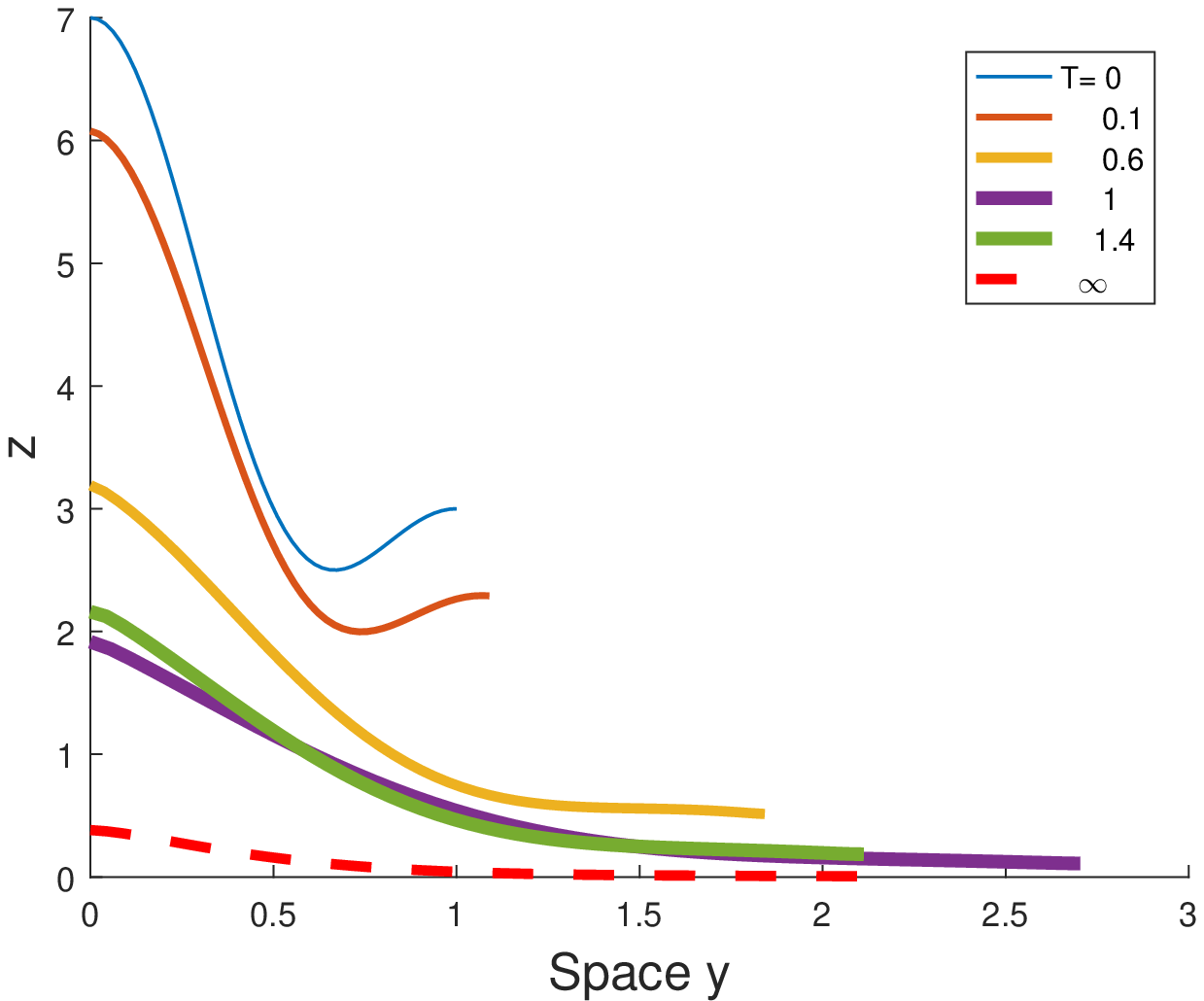}
} }
\subfigure[]{ {
\includegraphics[width=0.30\textwidth]{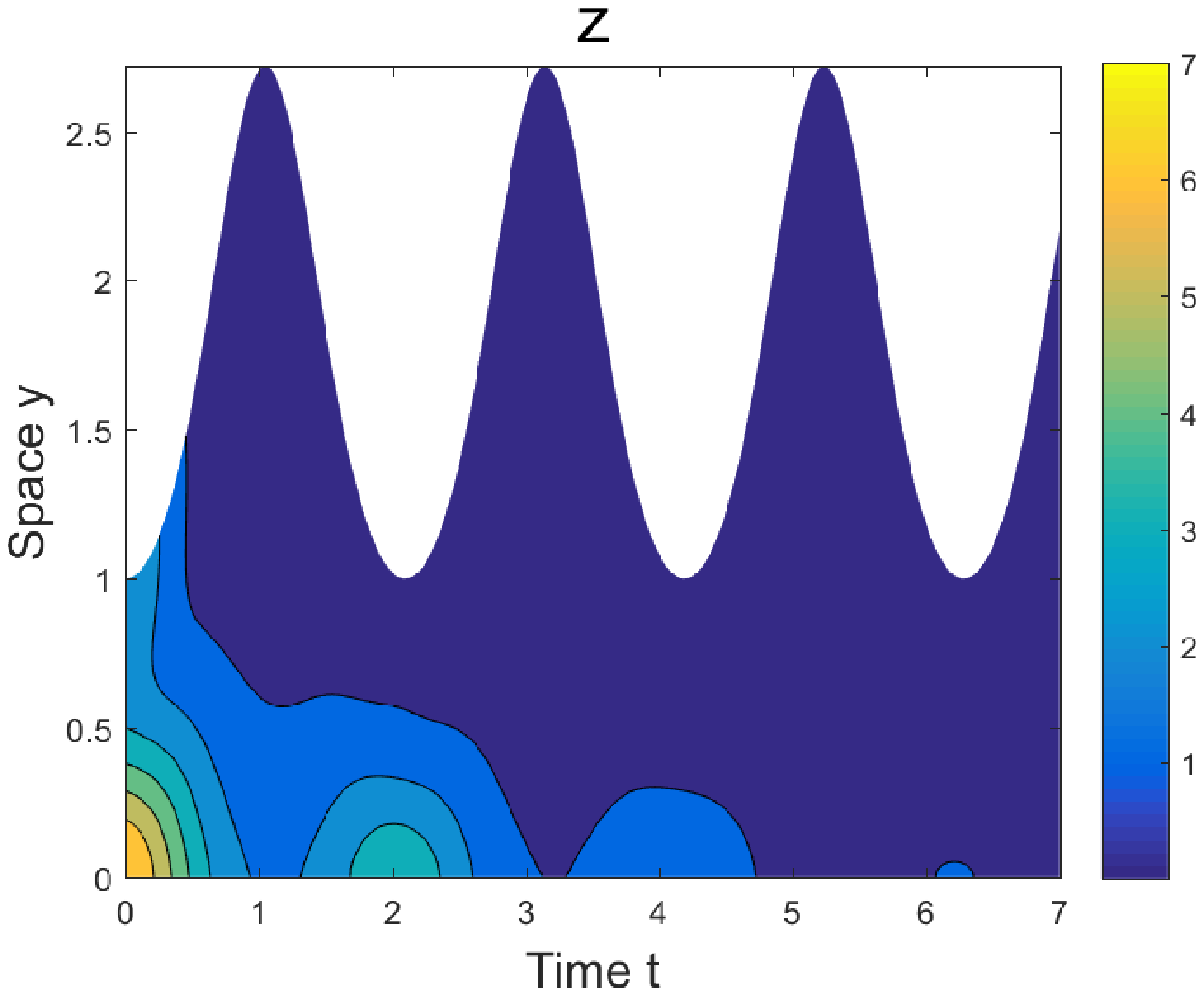}
} }
\caption{\scriptsize $\rho_4(t)=e^{-0.5(1-\cos(3t))}$. The interval is evolving with a larger evolution rate $\rho_4(t)$ and ${\mathcal{R}}_0<1$. Graphs $(a)$ show that $z$ is stable. Graphs $(b)$ and $(c)$, which are the cross-sectional view and contour map respectively, present the periodic evolution of the interval.}
\label{tu4}
\end{figure}

\begin{exm}
Fix $d(\rho(t)y,t)=0.1+0.2\rho(t)y$. We first choose $\rho_3(t)\equiv 1$, which means that the interval is fixed. Direct calculations show that
$${\mathcal{R}}_0^*(\rho_3)=\frac{\int_{0}^{T}\int_{0}^{L_0} g(I_0(t)e^{-k_0\rho(t)y})dydt}{\int_{0}^{T}\int_{0}^{L_0} d(\rho(t)y,t)dydt}
\approx\frac{0.4819}{0.4189}>1.$$
It is easily seen from Fig. \ref{tu3} that the phytoplankton species $z$ stabilizes to a positive periodic steady state.

Now we choose $\rho_4(t)=e^{0.5(1-\cos(3t))}$, it follows from \eqref{b31} that
$${\mathcal{R}}_0(\rho_4)=\frac{\int_{0}^{T}\int_{0}^{L_0} g(I_0(t)e^{-k_0\rho(t)y})dydt}{\int_{0}^{T}\int_{0}^{L_0} d(\rho(t)y,t)dydt}
\approx\frac{0.4536}{0.5767}<1.$$
Easily to see from Fig. \ref{tu4} that the phytoplankton species $z$ decays to zero.

\end{exm}

\begin{figure}[ht]
\centering
\subfigure[]{ {
\includegraphics[width=0.30\textwidth]{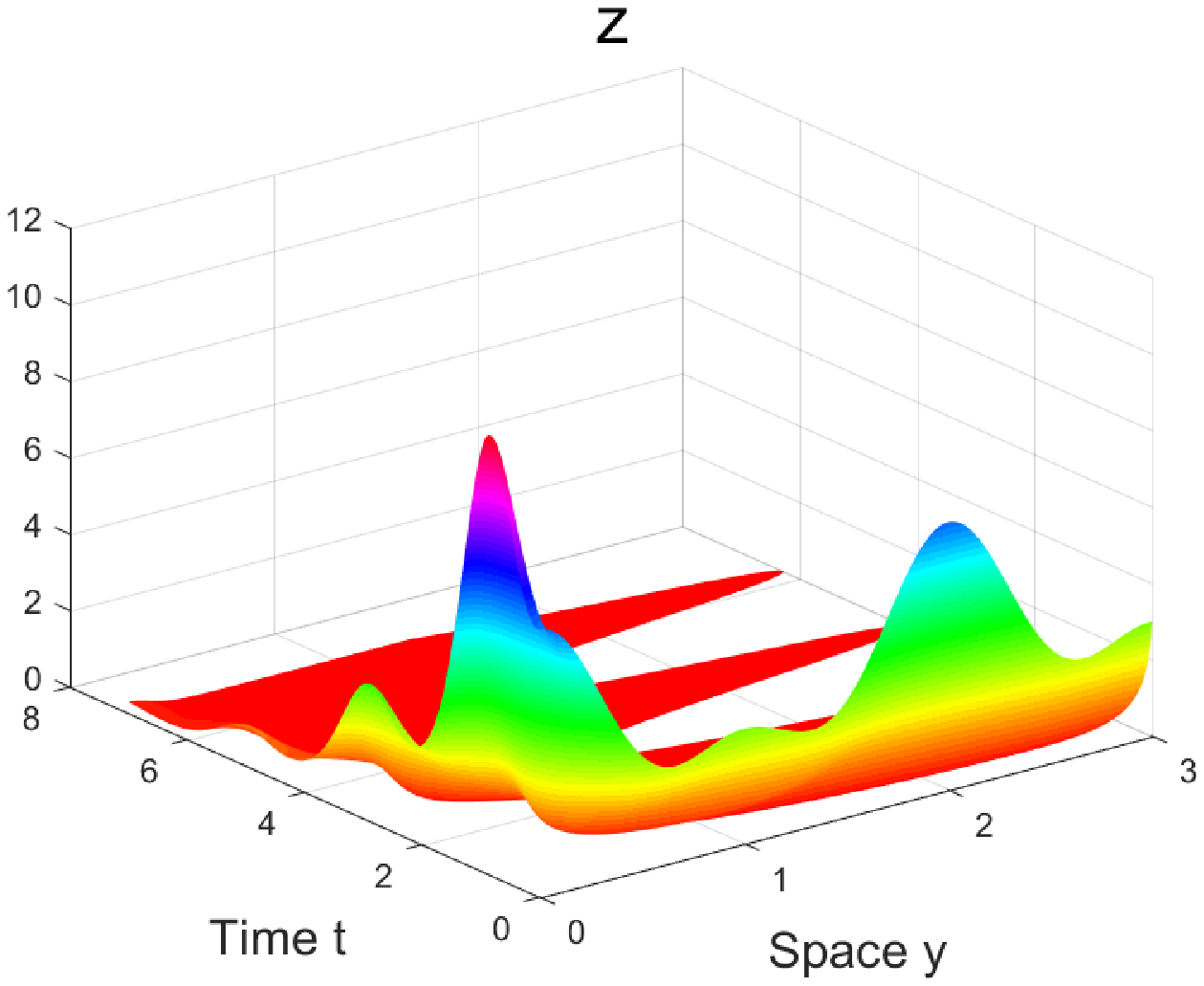}
} }
\subfigure[]{ {
\includegraphics[width=0.30\textwidth]{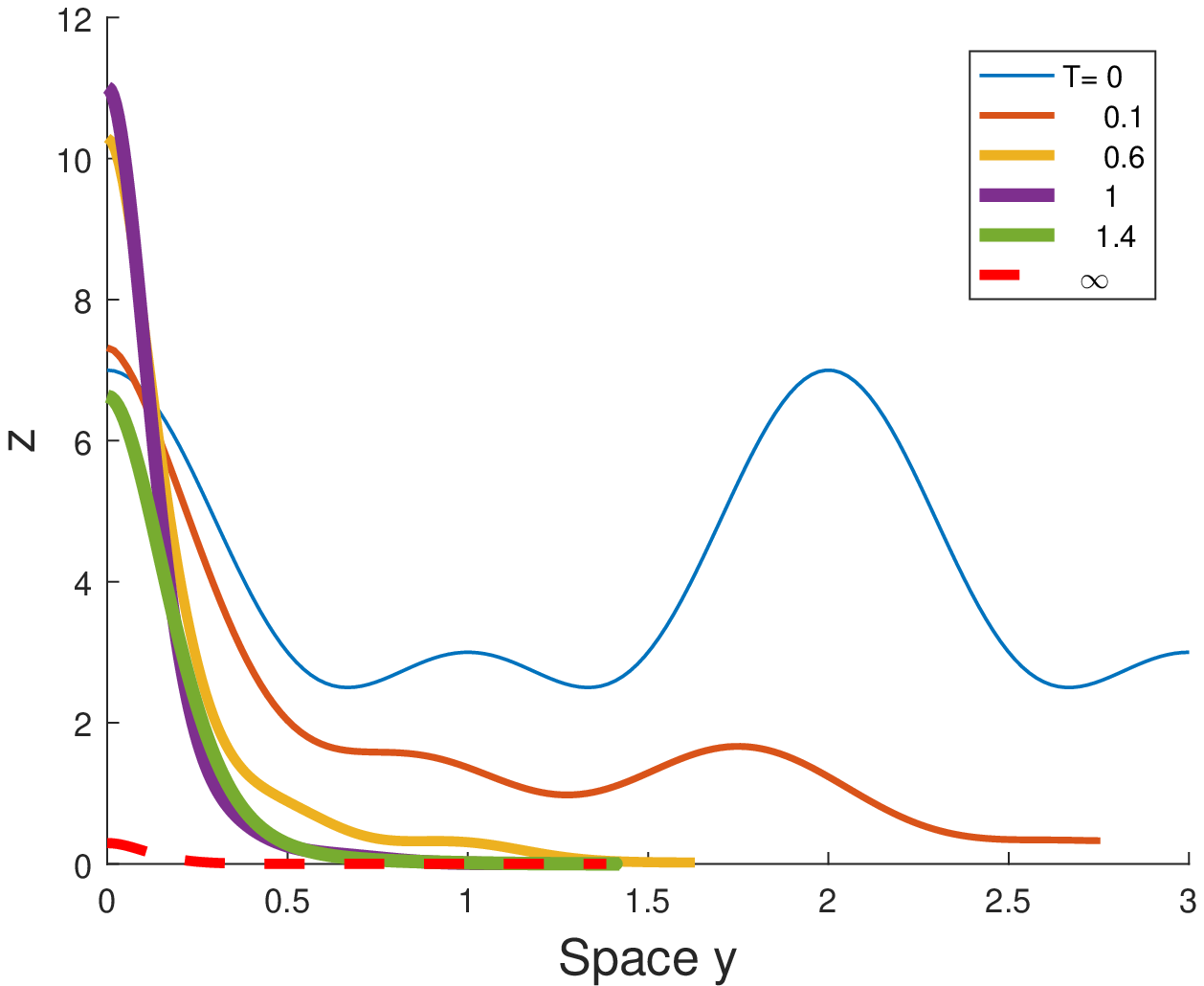}
} }
\subfigure[]{ {
\includegraphics[width=0.30\textwidth]{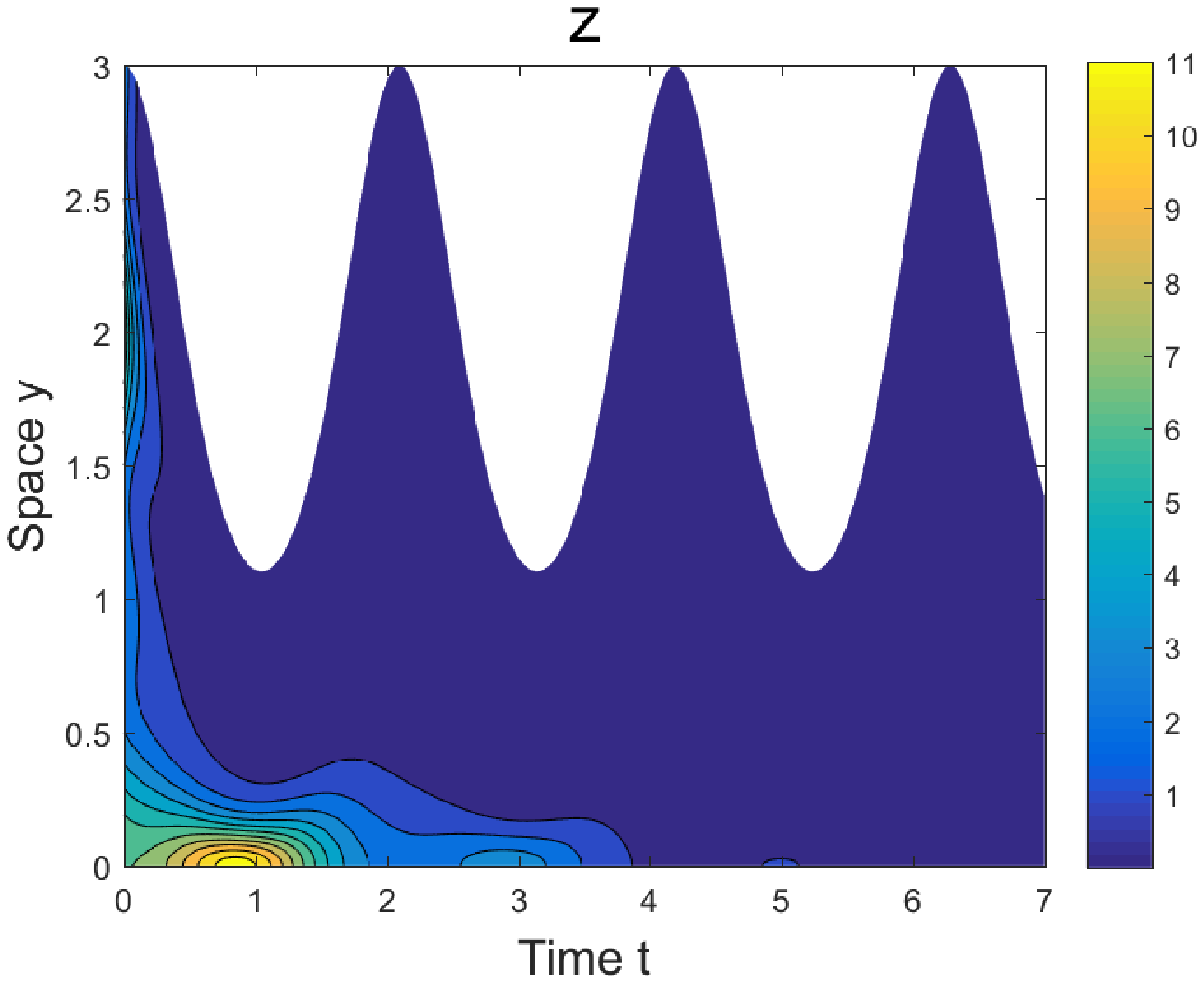}
} }
\caption{\scriptsize $d(\rho(t)y,t)=0.29+3\rho(t)y, L_0=3$ such that ${\mathcal{R}}_0<1$. Graphs $(a)-(c)$ showed that $z$ decays to $0$. }
\label{tu5}
\end{figure}

\begin{figure}[ht]
\centering
\subfigure[]{ {
\includegraphics[width=0.30\textwidth]{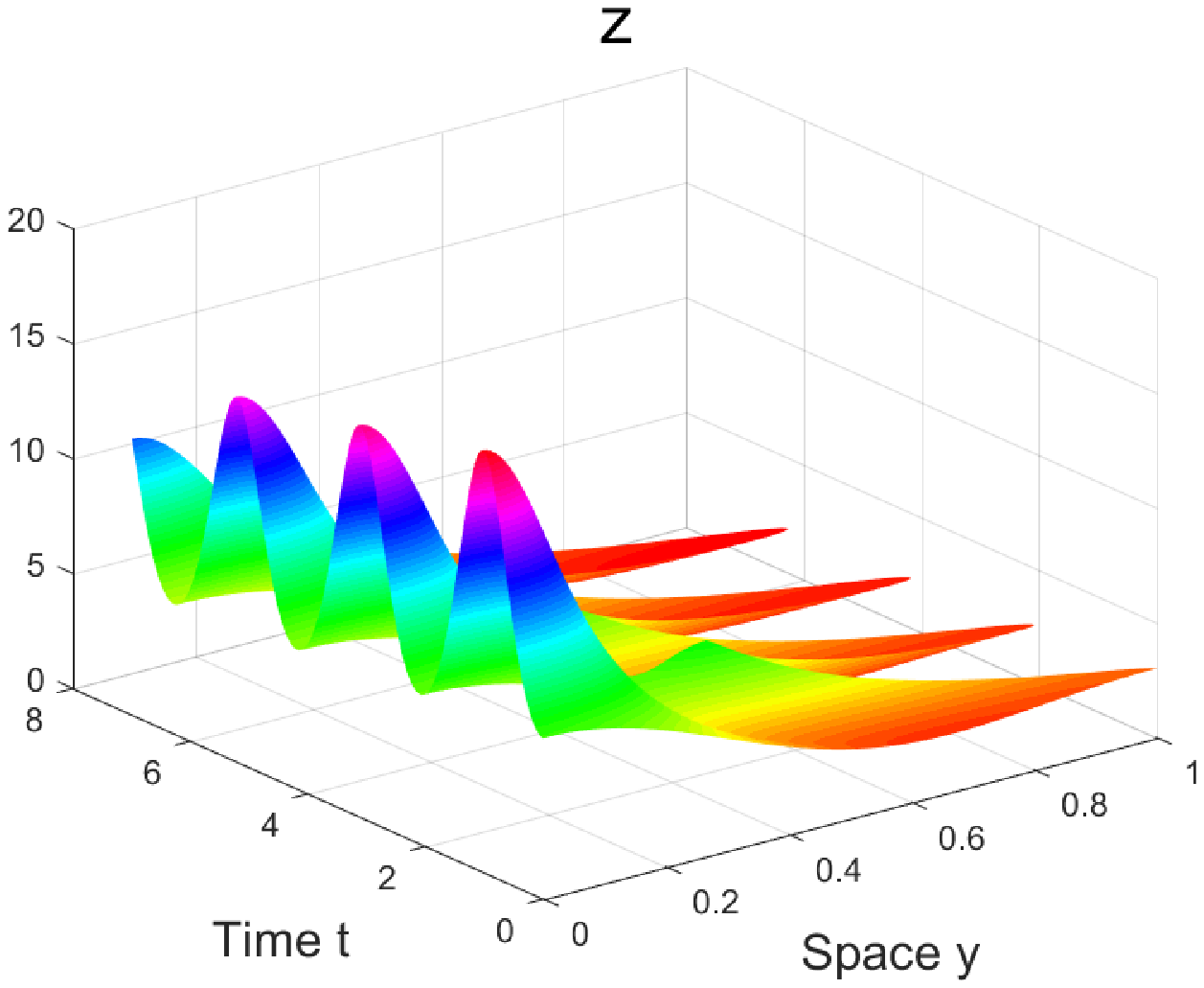}
} }
\subfigure[]{ {
\includegraphics[width=0.30\textwidth]{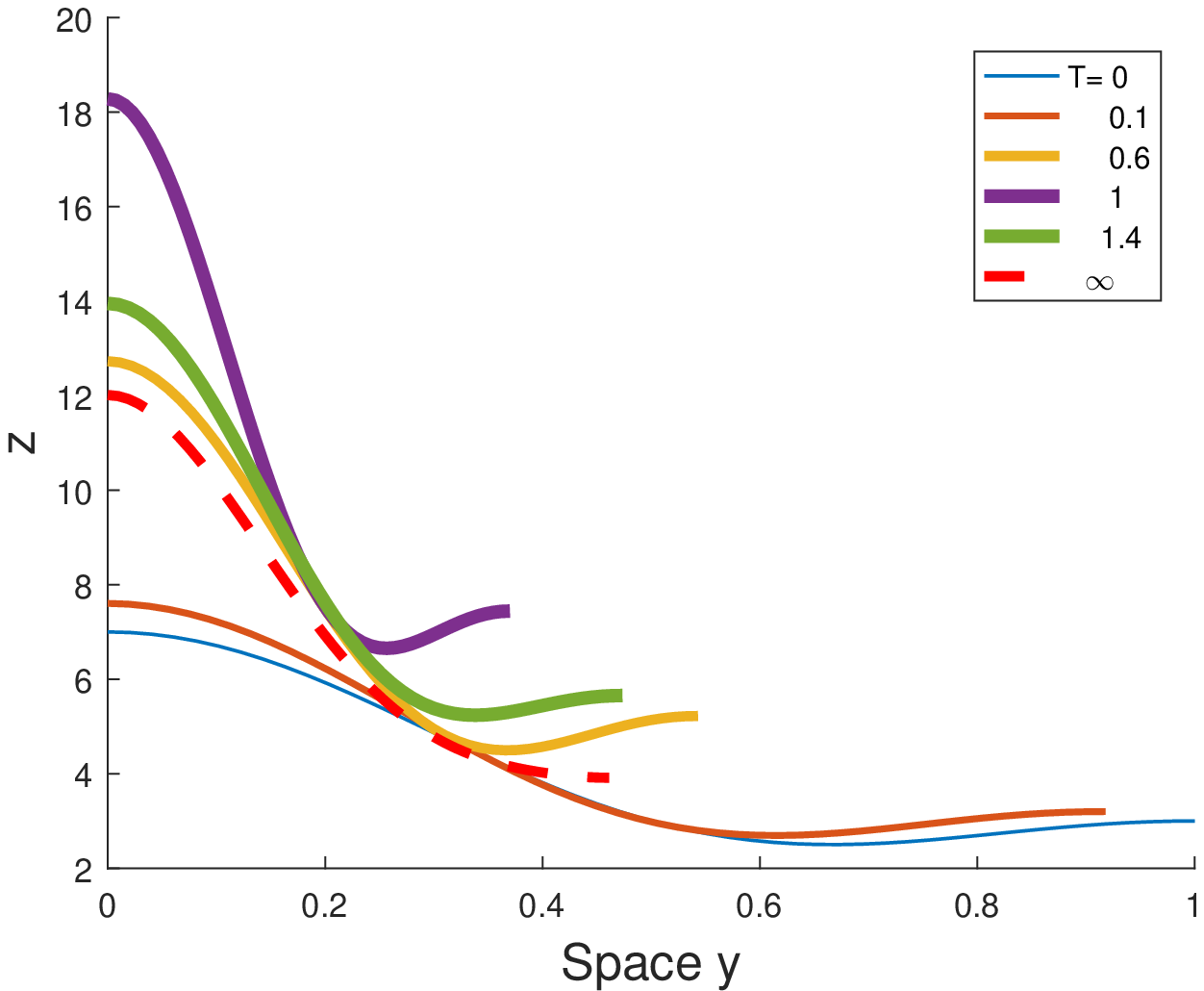}
} }
\subfigure[]{ {
\includegraphics[width=0.30\textwidth]{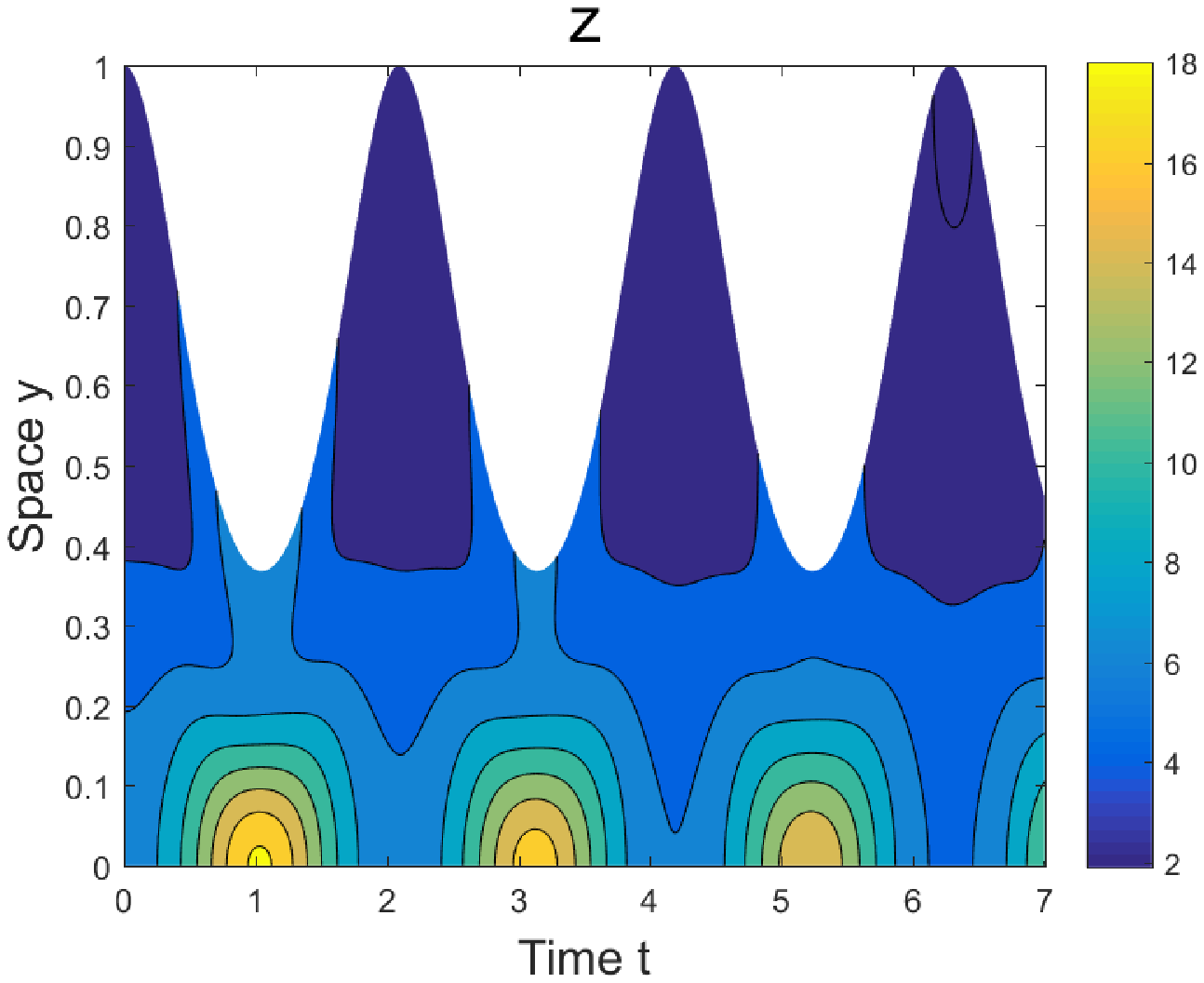}
} }
\caption{\scriptsize $d(\rho(t)y,t)=0.29+0.1\rho(t)y, L_0=1$ such that ${\mathcal{R}}_0>1$. Graphs $(a)-(c)$ showed that $z$ is stable. }
\label{tu6}
\end{figure}

Next, we verify the monotonicity of $\mathcal{R}_0$ with respect to $L$.
Assume that
\begin{equation*}
\left.
\begin{array}{lll}
D=0.001,\ \alpha=0, \  \ a=1,\ b=3,\   g=\frac{aI}{1+bI}=\frac{aI_0e^{-k_0\rho(t)y}}{1+bI_0e^{-k_0\rho(t)y}}, \\[7pt]
\ I_0=2,\ k_0=0.02, \ \rho_5(t)=e^{-0.5(1-\cos(3t))},\\[7pt]
\ z_0(y)=4+2\cos(\pi y)+\cos(2\pi y).
\end{array}
\right.
\end{equation*}
\begin{exm}
We first choose $d(\rho_5(t)y,t)=0.29+3\rho_5(t)y$ and $L_0=3$. From \eqref{b22}, direct calculations show that
$${\mathcal{R}}_0(\rho_5)\leq\frac{\int_{0}^{T}g(I_0(t))dt}{\int_{0}^{T}d(0,t)dt}
\approx\frac{0.5984}{0.6074}<1.$$
It is easily seen from Fig. \ref{tu5} that phytoplankton species $z$ decays to zero.

Now we choose $d(\rho_5(t)y,t)=0.29+0.1\rho_5(t)y$ and $L_0=1$, it follows from \eqref{b22} that
$${\mathcal{R}}_0(\rho_5)\geq\frac{\int_{0}^{T} g(I_0(t)e^{-k_0\rho(t)L_0})dt}{\int_{0}^{T} d(\rho(t)L_0,t)dt}
\approx\frac{0.5973}{0.579}>1.$$
Easily to get from Fig. \ref{tu6} that $z$ stabilizes to a positive periodic steady state.
\end{exm}

Recently, lots of models in evolving domains have been theoretically studied, for instance, models in period evolving domains \cite{SSM, wzc}, in growing domains \cite{tql1, tql2}, and in shrinking domains \cite{WN}.
In this paper, to understand the impacts of the evolving rate and depth on the dynamics of the single phytoplankton species, which depends upon the light for keeping the metabolism working constantly, a mathematical model in a periodically evolving environment has been established and studied.

The main conclusions are as follows: firstly, we introduce the basic reproducing number $\mathcal{R}_0$, which is used to describe threshold-type dynamics of our model, and give comparison principle in order to investigate the connection of $\mathcal{R}_0$ with regard to evolving rate $\rho(t)$, diffusion rate $D$, the sinking or buoyant rate $\alpha$ and the depth of water $L_0$. It is shown that that $\mathcal{R}_0(D,0,L)$ is strictly monotone decreasing in $L$ and $D$ respectively if $\alpha=0$. Especially, for enough small $D$ and $\alpha$, the limitation of $\mathcal{R}_0(D,\alpha,L)$ (denote $\mathcal{R}_0^*$) relies on the the domain evolution $\rho(t)$ and depth $L$. A point that should be stressed is that $\mathcal{R}_0^*$ is strictly monotone decreasing with respect to $\rho(t)$.
Secondly, it is proved in Theorem 3.1 that if  $\mathcal{R}_0\leq 1$, the phytoplankton species dies out in a long run, while if $\mathcal{R}_0>1$, then problem \eqref{a07}, \eqref{a09} admits at least one positive periodic solution, and there exists $\delta_0>0$ such that $\limsup_{t\to\infty}v(y,t)\geq\delta_0$ uniformly for $y\in [0,L_0]$, which means the phytoplankton species persists uniformly.

Furthermore, our numerical simulations indicate that the periodical domain evolution with small evolution rate has a positive effect on the survival of phytoplankton (see Figs. \ref{tu1} and \ref{tu2}), on the contrary, with large evolution rate
has a negative effect on the survival of phytoplankton (see Figs. \ref{tu3} and \ref{tu4}).
Meanwhile, large depth of water has a adverse effect on the survival of phytoplankton (see Figs. \ref{tu5} and \ref{tu6}), which matches closely our theoretical analysis.
In view of the complexity of the model, caused by $\rho(t)$ and $\alpha$, the mathematical expression of $\mathcal{R}_0$ has not yet given if $\alpha\neq0 $, we leave this for future work.

\end{document}